\documentclass[preprint,12pt]{elsarticle}

\usepackage{graphicx}
\usepackage{amssymb}
\usepackage{amsmath}
\usepackage{float}
\usepackage{tikz}
\usepackage{lineno}
\usetikzlibrary{positioning}
\usetikzlibrary{shapes,snakes}
\usetikzlibrary{decorations.markings}
\usepackage{freetikz}
\usepackage{booktabs}
\usepackage{subcaption}
\usepackage{diagbox}

\newcommand{\figref}[1]{Fig \ref{#1}}
\newcommand{\secref}[1]{Section \ref{#1}}
\newcommand{\tabref}[1]{Table \ref{#1}}
\newcommand{\appref}[1]{Appendix \ref{#1}}
\newcommand{\RNum}[1]{%
  \textup{\uppercase\expandafter{\romannumeral#1}}%
}
\renewcommand{\eqref}[1]{Eq.(\ref{#1})}
\DeclareMathOperator*{\argmin}{arg\,min}
\DeclareMathOperator*{\argmax}{arg\,max}

\journal{Reliability Engineering and System Safety}

\begin{document}

\begin{frontmatter}


\title{Optimal Inspection of Network Systems via Value of Information Analysis}

\author[1]{Chaochao Lin}
\ead{chaochal@andrew.cmu.edu}
\author[2]{Junho Song}
\ead{junhosong@snu.ac.kr}
\author[1]{Matteo Pozzi\corref{cor1}}
\ead{mpozzi@cmu.edu}
\address[1]{Dept. of Civ. and Env. Eng., Carnegie Mellon University, Pittsburgh, PA, USA}
\address[2]{Dept. of Civ. and Env. Eng., Seoul National University, Seoul, South Korea}
\cortext[cor1]{Corresponding author}

% \linenumbers

\begin{abstract}
This paper develops computable metrics to assign priorities for information collection on network systems made up by binary components. Components are worth inspecting because their condition state is uncertain and the system functioning depends on it. The Value of Information (VoI) allows assessing the impact of information in decision making under uncertainty, including the precision of the observation, the available actions and the expected economic loss. Some VoI-based metrics for system-level and component-level maintenance actions, defined as “global” and “local” metrics, respectively, are introduced, analyzed and applied to series and parallel systems. Their computationally complexity of applications to general networks is discussed and, to tame the complexity for the local metric assessment, a heuristic is presented and its performance is compared on some case studies.
\end{abstract}


\end{frontmatter}

%\begin{linenumbers}

\section{Introduction}
\label{sec:introduction}
Many civil infrastructures consist of multiple binary components, arranged as a network to fulfill the function of the system. The binary states of the components, either intact or damaged, determine the system condition. 
The belief of the agent controlling the maintenance process can be described by a probabilistic distribution on the possible states of the components. Maintenance actions are selected to trade off the risk of system malfunctioning with the cost of maintenance (repair and retrofitting).
Observations of the components' states can improve decision making and reduce the uncertainty and maintenance cost. However, due to budget constraints, it is often impossible to inspect all components in the network. Therefore it is important to assign inspection priorities among components. Intuitively, many factors can affect the inspection preferences, such as the probabilities of failure events, the maintenance costs and the role of each component in the system's functions. These factors can be integrated in an Importance Measure (IM), a value assigned to each component to summarize the benefit of inspection.

To introduce the problem, consider a binary system made up by $N$ components: $\{c_1,c_2,\cdots,c_N\}$.
Let $\mathit{s}=[\mathit{s}_1, \mathit{s}_2, \cdots, \mathit{s}_N]\in S$ denote the states of the components, with $\mathit{s}_j=1$ indicating that component $c_j$ is working, and $s_i=0$ that it fails, where $S=\mathbb{B}^N$ and $\mathbb{B}=\{0,1\}$. The system state $u=\phi(\mathit{s})$ is also a binary variable, where $\phi:S\rightarrow \mathbb{B}$ is the component-to-system function. State $\mathit{s}$ is unknown to the agent who manages the system. Instead, the agent optimizes the measurement and maintenance plans based on her belief of $\mathit{s}$. The prior probability distribution of $\mathit{s}$ is denoted as $p_\mathit{s}:S \rightarrow [0,1]$, and $p_j=\mathbb{P}[\mathit{s}_j=0]$ indicates the prior marginal failure probability of $c_j$. The failure probability of the system is $p_u=\mathbb{P}[u=0]$, and we use $p_\pi$ and $p_{\omega|E}$ for the prior value of $p_u$ and its posterior value given event $E$, respectively.

In this paper, we develop metrics to assess the importance of inspecting any component. We assume that the outcome of the inspection is also binary. If component $c_i$ is inspected, $y_i=0$ indicates an “alarm”, i.e. a symptom that $c_i$ is not working, while $y_i=1$ indicates that $c_i$ seems to work, and we define this outcome as a “silence”. If the inspection is perfect, then $y_i=s_i$.

Based on the measurement outcome, the prior distribution of random variables $\mathit{s}$ can be updated to posterior distribution $p_{\mathit{s}|\mathit{y}_i}$ and, the system level failure probability to $p_{\omega|\mathit{y}_i}$. When the components are interdependent, the measurement on one component may also affect the failure probability on other components.

Importance Measures (IMs) were first introduced by Birnbaum \citep{birnbaum1968importance} to evaluate the contribution of each component to the system performances, as network connectivity. The Birnbaum's Measure (BM) \citep{birnbaum1968importance} evaluates the importance of a component by the difference in the posterior system failure probability when it is damaged or intact (i.e., in our framework, when the inspection outcome is alarm or silence):
\begin{equation}
\text{BM}(i)= p_{\omega|\mathit{y}_i=0} - p_{\omega|\mathit{y}_i=1}
\end{equation}
Other IMs are briefly discussed in the \appref{app:Importance_Measures}. Most of them focus on the marginal or conditional probability of the failure events, and do not explicitly include any evaluation of the maintenance cost and risk related to the problem.  In maintenance problems, some components need high attention because of their topological role in the network, some others because of the high risk that the component fails. To address this issue, \citep{Wu2013} extends the BM to a new cost-based IM. \cite{zio2007importance} provides an approach which achieves multi-objective (such as system risk and economics) optimization problem and generic algorithms to reduce the computation time. \cite{DerKiureghian2007} modeled the component failures as independent Poisson events and developed IMs for long-term maintenance of series, parallel and general systems based on the system unavailability, mean rate of failure and mean duration of downtime.

To compare and rank the impact of inspections, one can assess their Value of Information (VoI). VoI assessment is based on Bayesian pre-posterior analysis, as introduced by \citep{Howard1966InformationValue}, integrating the probabilistic knowledge about the system with the economic factors related to the available actions. In the maintenance process of infrastructure systems, the economic costs are related to the system malfunctioning, the execution of inspections, repairing or replacement actions, etc. VoI has been studied intensively in the area of Structural Health Monitoring (SHM). \citep{Pozzi2011AssessingMonitoring} provides a general framework for assessing VoI for the long-term SHM, and proposes a Monte Carlo approach to reduce the computation complexity. \cite{Pozzi2011AssessingMonitoring} also investigates how the imperfect measurements affect the posterior decisions.
\citep{STRAUB2005335} integrates VoI with risk based inspection to schedule inspections and maintenance planning of structural systems.
\citep{straub2017value} investigates the stochastic dependencies in component deterioration, failure consequence of the system state, the component inspection cost and performance in structural systems and how they will affect the VoI distribution.

VoI has also been applied to long-term decision making problems. \cite{MillerValueSeqInfo} extends the VoI analysis to optimize not only static one-shot inspection, but also sequentially dependent observations. \cite{memarzadeh2016value} and \cite{andriotis2019value} apply the component-wise VoI metric to sequential decision making in the management of infrastructure systems, modeled by Partially Observable Markov Decision Process (POMDP). \cite{thons2018value} uses decision trees to assess long-term VoI.

The complexity of computing VoI can grow exponentially as the number of components in the system increases \citep{malings2018value}. Even worse, the VoI generally lacks the property of submodularity \cite{Malings2019}, so that the application of greedy approaches does not provide certain guarantees of near-optimal solutions \citep{shamaiah2010greedy}. Simplifications have been proposed for efficient VoI computation in some special cases  \citep{MALINGS201677}. 

In this paper, we investigate VoI-based metrics related to system-level (“global”) and component-level (“local”) decision making after component inspections, for networks with various topologies, and compare these results with traditional IMs.
A recent paper, \cite{Fauriat2019}, also focuses on inspections for networked systems, and it develops approaches to identify the components most in need of inspections, adopting an approach similar to what we define the local metric.
We also derive simple optimal rules for series and parallel systems. For general networks, the computational complexity of the problem is discussed and a heuristic approach is provided. \secref{sec:voi} introduces the global and local metrics for evaluating the components' VoI. \secref{sec:optimal_inspection} identifies optimal rules for optimizing these metrics to typical networks such as series and parallel systems. \secref{sec:heuristic} proposes an approximated approaches to simplify the optimization complexity, and \secref{sec:examples} explores different applications of global, local and heuristic approaches to some network examples.
\section{Global and local VoI metrics}
\label{sec:voi}

\subsection{Principles of VoI}
\label{sub:general_voi}
\begin{figure}
    \centering
    \begin{subfigure}[b]{0.4\textwidth}
    \begin{tikzpicture}[
        roundnode/.style={circle, draw=black!100, very thick, minimum width=7mm},
        diamondnode/.style={diamond, draw=black!100, very thick, minimum width=5mm},
        squarenode/.style={rectangle, draw=black!100, very thick, minimum size=7mm},
        ]
        %Nodes
        \node[roundnode] (maintopic) {$\mathit{s}$};
        \node[roundnode] (nextstate) [right=of maintopic] {$\mathit{s'}$};
        \node[roundnode] (measure) [below= of maintopic] {$\mathit{y_i}$};
        \node[squarenode] (action) [right=of measure] {$A$};
        \node[roundnode] (up) [right=of nextstate] {$u'$};
        \node[diamondnode] (loss) [right=of up] {$\mathcal{L}_\RNum{1}$};
        \node[diamondnode] (lossp) [right=27mm of action] {$\mathcal{L}_\RNum{2}$};
         
        %Lines
        \draw[
            decoration={markings,mark=at position 1 with {\arrow[scale=2]{>}}},
            postaction={decorate},
            shorten >=0.4pt
            ] (maintopic.east) -- (nextstate.west);
        \draw[
            decoration={markings,mark=at position 1 with {\arrow[scale=2]{>}}},
            postaction={decorate},
            shorten >=0.4pt
            ] (maintopic.south) -- (measure.north);
        \draw[
            decoration={markings,mark=at position 1 with {\arrow[scale=2]{>}}},
            postaction={decorate},
            shorten >=0.4pt, dashed
            ] (measure.east) -- (action.west);
        \draw[
            decoration={markings,mark=at position 1 with {\arrow[scale=1]{>}}},
            postaction={decorate},
            shorten >=0.4pt, double
            ] (nextstate.east) -- node[anchor=south] {$\phi$} (up.west);
        \draw[
            decoration={markings,mark=at position 1 with {\arrow[scale=2]{>}}},
            postaction={decorate},
            shorten >=0.4pt
            ] (up.east) -- (loss.west);
        \draw[
            decoration={markings,mark=at position 1 with {\arrow[scale=2]{>}}},
            postaction={decorate},
            shorten >=0.4pt
            ] (action.east) -- (lossp.west);
        \draw[
            decoration={markings,mark=at position 1 with {\arrow[scale=2]{>}}},
            postaction={decorate},
            shorten >=0.4pt
            ] (action.north) -- (nextstate.south);
    \end{tikzpicture}
    \caption{}
    \label{fig:pgma}
    \end{subfigure}
    \hfill
    \begin{subfigure}[b]{.4\textwidth}
        \begin{tikzpicture}[
            roundnode/.style={circle, draw=black!100, very thick, minimum width=7mm},
            diamondnode/.style={diamond, draw=black!100, very thick, minimum width=5mm},
            squarenode/.style={rectangle, draw=black!100, very thick, minimum size=7mm},
            ]
            %Nodes
            \node[roundnode] (maintopic) {$\mathit{u}$};
            \node[roundnode] (nextstate) [right=of maintopic] {$\mathit{u'}$};
            \node[roundnode] (measure) [below= of maintopic] {$\mathit{y_i}$};
            \node[squarenode] (action) [right=of measure] {$A$};
            \node[diamondnode] (loss) [right=of nextstate] {$\mathcal{L}_\RNum{1}$};
            \node[diamondnode] (lossp) [right=of action] {$\mathcal{L}_\RNum{2}$};
             
            %Lines
            \draw[
                decoration={markings,mark=at position 1 with {\arrow[scale=2]{>}}},
                postaction={decorate},
                shorten >=0.4pt
                ] (maintopic.east) -- (nextstate.west);
            \draw[
                decoration={markings,mark=at position 1 with {\arrow[scale=2]{>}}},
                postaction={decorate},
                shorten >=0.4pt
                ] (maintopic.south) -- (measure.north);
            \draw[
                decoration={markings,mark=at position 1 with {\arrow[scale=2]{>}}},
                postaction={decorate},
                shorten >=0.4pt, dashed
                ] (measure.east) -- (action.west);
            \draw[
                decoration={markings,mark=at position 1 with {\arrow[scale=2]{>}}},
                postaction={decorate},
                shorten >=0.4pt
                ] (nextstate.east) -- (loss.west);
            \draw[
                decoration={markings,mark=at position 1 with {\arrow[scale=2]{>}}},
                postaction={decorate},
                shorten >=0.4pt
                ] (action.east) -- (lossp.west);
            \draw[
                decoration={markings,mark=at position 1 with {\arrow[scale=2]{>}}},
                postaction={decorate},
                shorten >=0.4pt
                ] (action.north) -- (nextstate.south);
        \end{tikzpicture}
        \caption{}
    \label{fig:pgmb} 
    \end{subfigure}
    \caption{Decision graph for the general problem (a), and for the global metric (b).}
    \label{fig:figure1}
\end{figure}
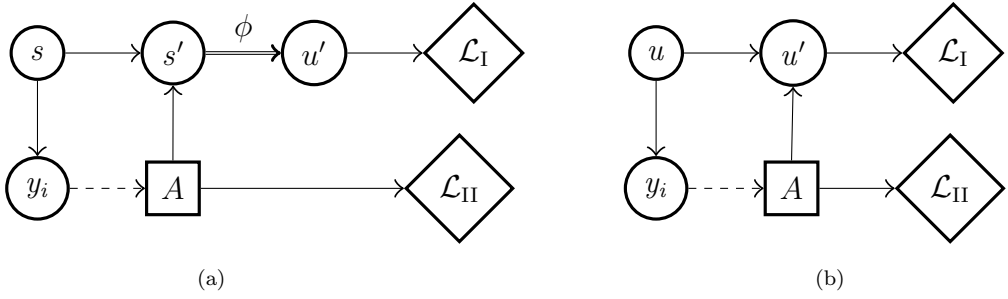

The decision graph for the process of inspecting and maintaining the system is illustrated in \figref{fig:pgma}.
Let $\mathcal{A}$ denote the set of all possible maintenance plans, that we simply call “actions”. Action $A\in \mathcal{A}$ transforms current components' state $\mathit{s}\in S$ into state $\mathit{s'}\in S$, via transition distribution $p_{\mathit{s'}|\mathit{s},A}:S\times \mathcal{A}\times S\rightarrow [0,1] $.
Loss function $\mathcal{L}(\mathit{s'}, A)=\mathcal{L}_\RNum{1}(\phi(\mathit{s'}))+\mathcal{L}_\RNum{2}(A):S\times \mathcal{A}\rightarrow \mathbb{R}$ summarizes the overall cost: $\mathcal{L}_\RNum{1}(\phi(\mathit{s'}))=C_F(1-u')$ adds failure costs $C_F$ if the system is not functioning, depending on system state $u'$ after taking action $A$, which is associated with implementing cost $\mathcal{L}_\RNum{2}(A)$.

The prior loss $L_\pi$ is the minimum expected cost among all possible actions, before any inspection:
\begin{equation}
    \label{eq:prior_loss}
    L_\pi = \min_{A}\mathbb{E}_\mathit{s}\mathbb{E}_{\mathit{s'}|\mathit{s},A} \mathcal{L}(\mathit{s'},A)=\min_{A} \mathbb{E}_{\mathit{s'}|A} \mathcal{L}(\mathit{s'},A)
\end{equation}
where $\mathbb{E}_{\mathit{s'}|A}[\cdot]=\mathbb{E}_s\mathbb{E}_{s'|s,A}[\cdot]$ denotes the statistical expectation depending on distributions $p_{\mathit{s'}|\mathit{s},A}$ and $p_{\mathit{s}}$.

Inspecting component $c_i$, the agent collects observation $y_i$ distributed according to function $p_{\mathit{y}_i}:\mathbb{B}\rightarrow [0,1]$, and the belief of the components' state $\mathit{s}$ is updated to posterior distribution $p_{\mathit{s}|\mathit{y}_i}:S\times \mathbb{B}\rightarrow [0,1]$. These functions are obtained by Bayes' rule:
\begin{equation}
    \label{eq:Bayes}
    p_{\mathit{y}_i}= \sum_{\mathit{s}}  p_{\mathit{y}_i|\mathit{s}} p_{\mathit{s}}\;\;\;\;
    p_{\mathit{s}|\mathit{y}_i}=\frac{ p_{\mathit{y}_i|\mathit{s}} p_{\mathit{s}} } { p_{\mathit{y}_i} }
\end{equation}
where $p_{\mathit{y}_i|\mathit{s}}:\mathbb{B}\times S\rightarrow [0,1]$ is the likelihood function related to observation $\mathit{y}_i$.

The corresponding expected posterior loss is:
\begin{equation}
    \label{eq:posterior_loss}
    L_\omega(i)=\mathbb{E}_\mathit{y_i}\min_A\mathbb{E}_{\mathit{s'}|\mathit{y}_i,A}\mathcal{L}(\mathit{s'},A)
\end{equation}
where $\mathbb{E}_{\mathit{s'}|\mathit{y}_i,A}[\cdot]=\mathbb{E}_{\mathit{s}|\mathit{y}_i}\mathbb{E}_{\mathit{s'}|\mathit{s},A}[\cdot]$ is the posterior expectation, related to distribution $p_{\mathit{s}|\mathit{y}_i}$, and $\mathbb{E}_\mathit{y_i}[\cdot]$ is  related to distribution $p_{\mathit{y}_i}$.

The VoI for inspecting $c_i$ is the expected loss reduction due to the inspection, i.e. the difference of the prior and posterior loss functions \cite{Howard1966InformationValue}:
\begin{equation}
    \label{eq:voiy}
    \text{VoI}(i)=L_\pi - L_\omega(i)
\end{equation}

Loss function $\mathcal{L}$ does not include the cost of monitoring, and the VoI is always not negative. However, if such cost is uniform among components, the VoI is a rational IM assessing the relevance of inspections. 
The optimal component to inspect, $c_{i^\star}$, is the argument that maximizes \eqref{eq:voiy}:
\begin{equation}
    \label{eq:ystar}
    \mathit{i}^\star = \argmax_\mathit{i} \text{VoI}(\mathit{i})
\end{equation}

The VoI depends on the specific number $N$ of components, the action domain $\mathcal{A}$, the loss function $\mathcal{L}$ (in turn defined by the component-to-system function $\phi$, the failure cost $C_F$, and the implementing cost $\mathcal{L}_\RNum{2}$), the prior probability $p_{\mathit{s}}$, the transition probability $p_{\mathit{s'}|\mathit{s},A}$ and the likelihood function $p_{\mathit{y}_i|\mathit{s}}$ adopted, as apparent in \figref{fig:pgma}. In the following Sections, we describe a form of the likelihood function for binary components, and then we focus on two classes of losses and transitions, related to global and local decision making.

\subsection{Modeling imperfect inspections}
\label{sec:measurement_error}

The VoI analysis also depends on the specific assumed likelihood function. If the binary outcome $\mathit{y}_i$, of inspecting component $c_i$, only depends on the state $\mathit{s}_i$ of that component, likelihood function $p_{\mathit{y}_i|\mathit{s}}$ in \eqref{eq:Bayes} is reduced to a 4-entry emission table $p_{\mathit{y}_i|\mathit{s}_i}:\mathbb{B}\times\mathbb{B}\rightarrow[0,1]$, shown in \tabref{tab:false_rate}.

Observations of components' states are prone to error, and the inaccuracy can be formulated by two parameters, $\epsilon_{\text{FS}}=\mathbb{P}[\mathit{y}_i=1|\mathit{s}_i=0]$ and $\epsilon_{\text{FA}}= \mathbb{P}[\mathit{y}_j=0|\mathit{s}_j=1]$, which are the probability of type $\RNum{1}$ error: having a “silence” when the component is undamaged, and of type $\RNum{2}$ error: an alarm when the component is damaged, respectively. While these probabilities can depend on the specific component, in the following we assume that all the components have identical $\epsilon_{\text{FS}}$ and $\epsilon_{\text{FA}}$.

\begin{table}
    \centering
    \begin{tabular}{lcc}\hline
        \toprule
        \diagbox[width=15em]{Actual state}{Observation}&
          Silence $y_i=1$ & Alarm $y_i=0$ \\
        \midrule
        Undamaged $s_i=1$ & $1-\epsilon_{\text{FA}}$ & $\epsilon_{\text{FA}}$ \\
        \midrule
        Damaged $s_i=0$ & $\epsilon_{\text{FS}}$  &  $1-\epsilon_{\text{FS}}$ \\
        \bottomrule
        \end{tabular}
        \caption{Emission probability table for observation $\mathit{y}_i$ given state $\mathit{s}_i$.}
        \label{tab:false_rate}
\end{table}

Inspection outcomes probability function $p_{\mathit{y}_i}:\mathbb{B}\rightarrow[0,1]$, is related to a single value: the probability $h_i=\mathbb{P}[\mathit{y}_i=0]$ of receiving an alarm on $c_i$, which is:
\begin{equation}
    \label{eq:probAlarm}
        h_i =(1-\epsilon_{\text{FS}})p_i + \epsilon_{\text{FA}}(1-p_i) = \epsilon_{\text{FA}} + K p_i
\end{equation}
where constant $K=1-\epsilon_{\text{FA}}-\epsilon_{\text{FS}}$ is strictly positive, as we assume that both $\epsilon_{\text{FA}}$ and $\epsilon_{\text{FS}}$ are less than $1/2$.

\subsection{Global metric}
\label{sub:global_metric}
\begin{figure}
    \centering
    \begin{subfigure}[b]{0.45\textwidth}
        \includegraphics[width=\textwidth]{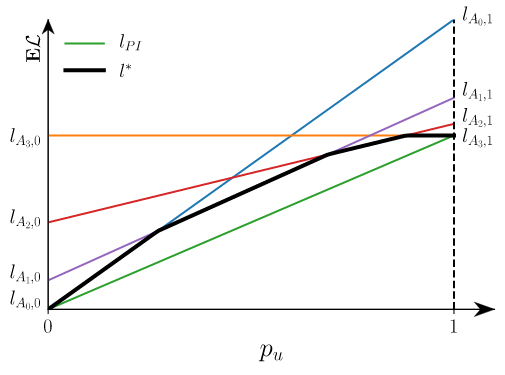}
        \caption{}
        \label{fig:global_general_a}
    \end{subfigure}
    \hfill
    \begin{subfigure}[b]{0.45\textwidth}
        \includegraphics[width=\textwidth]{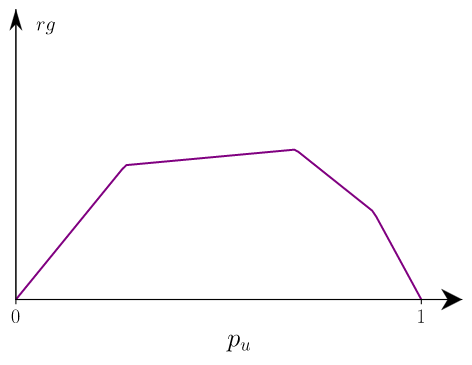}
        \caption{}
        \label{fig:global_general_b}
    \end{subfigure}
    \caption{Expected loss function (a) and corresponding regret (b) for a global problem with 4 possible actions.}
    \label{fig:global_general}
\end{figure}

We define the global metric assuming that action $A$ affects the system state $u$. In this setting, for any of the two values of the binary variable $u$, an expected loss value can be assigned to any action $A$, regardless of the details of components' condition described by variable $\mathit{s}$. \figref{fig:pgmb} shows the corresponding decision graph, where the loss is a function of system state $u'$ after the taken action: $l(u',A)=\mathcal{L}(s',A)$, with $u'=\phi(s')$. Transition function $p_{\mathit{s'}|\mathit{s},A}$ is now converted into function $p_{u'|u,A}:\mathbb{B}\times\mathcal{A}\times\mathbb{B}\rightarrow[0,1]$, in turn defined by the pair of values, $p_{0,A,0}$ and $p_{1,A,0}$, which are the probabilities that $u'=0$ given action $A$ and given $u=0$ or $u=1$, respectively. Then, $l_{A,0}=p_{0,A,0} C_F+C_A$ and $l_{A,1}=p_{1,A,0} C_F+C_A$, with $C_A=\mathcal{L}_\RNum{2} (A)$ represents the expected loss when $u$ is zero and one (i.e. when the system is not working and is working), respectively, for action $A$.
For each pair of values $p_{1,A,0}\leq p_{0,A,0}$ so that $l_{A,0}-l_{A,1}\leq C_F$, one can find a pair of values $C_A=l_{A,1}$ and $p_{0,A,0}=(l_{A,1}-l_{A,0})/C_F$, to represent the target pair of losses, assuming that no maintenance action makes the system degrade, so $p_{1,A,0}=0$. 
The agent has to find an optimal trade-off between implementing more expensive actions related to a low risk $p_{0,A,0}$, and less expensive ones, related to a higher risk.

The corresponding expected loss under action $A$ is a linear function of the system failure probability $p_u$:
\begin{equation}
    l_A (p_u )=\mathbb{E}_u\mathbb{E}_{u'|u,A} l(u',A)=p_u l_{A,0}+(1-p_u ) l_{A,1}
\end{equation}
By taking the minimum among available actions in domain $\mathcal{A}$, the optimal loss is defined by concave function $l^* (p_u )=\min_A l_A (p_u )$. The prior expected loss of \eqref{eq:prior_loss} for the global metric is thus $L_\pi^\text{G}=l^* (p_\pi )$ and, following \eqref{eq:prior_loss}, the posterior loss inspecting $c_i$ is:
\begin{equation}
    \label{eq:posterior_loss_global}
    L_\omega^\text{G} (i)=h_i l^* (p_{\omega|y_i=0})+(1-h_i ) l^* (p_{\omega|y_i=1})
\end{equation}
and the VoI, following \eqref{eq:voiy}, is $\text{VoI}_\text{G} (i)=L_\pi^\text{G}-L_\omega^\text{G} (i)$.

\begin{figure}
    \centering
    \begin{subfigure}[b]{0.45\textwidth}
        \includegraphics[width=\textwidth]{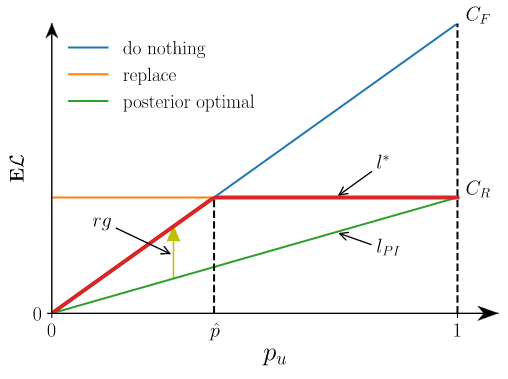}
        \caption{}
        \label{fig:global_one}
    \end{subfigure}
    \hfill
    \begin{subfigure}[b]{0.45\textwidth}
        \includegraphics[width=\textwidth]{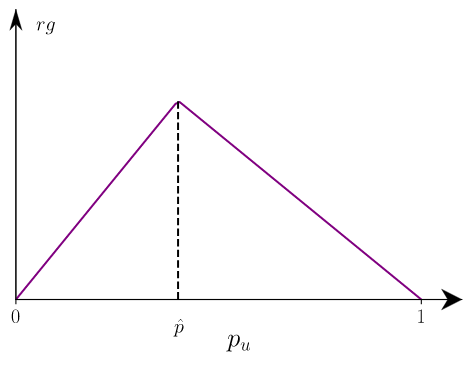}
        \caption{}
        \label{fig:rg}
    \end{subfigure}
    \caption{Expected loss (a) and corresponding regret (b) for the binary actions case.}
    \label{fig:expected_loss_rg}
\end{figure}

As a function of $p_u$, the expected loss with perfect information of $u$ is linear function $l_{\text{PI}} (p_u)=p_u l_0^*+(1-p_u) l_1^*$, with $l_0^*=\min_A l_{A,0}=l^*(1)$ and $l_1^*=\min_A l_{A,1}=l^* (0)$, and the “regret” is the concave function $rg(p_u )=l^*(p_u)-l_{\text{PI}} (p_u)$, with $rg(0)=rg(1)=0$. 
The corresponding prior regret is $\text{RG}_\pi=rg(p_\pi )$.
As function $l_{\text{PI}}$ is linear, the expected posterior loss with perfect information is $L_{\text{PI}}=l_{\text{PI}} (p_\pi)$, and expected posterior regret inspecting $c_i$ is $\text{RG}_\omega(i)=L_\omega^\text{G} (i)-L_{\text{PI}}=-\text{VoI}_\text{G} (i)+L_\pi^\text{G}-L_{\text{PI}}$.
Hence, component $c_{i^*}$, that maximizes the VoI identified in \eqref{eq:ystar}, also minimizes the expected posterior regret:
\begin{equation}
    \mathit{i}^\star = \argmin_\mathit{i} \text{RG}_\omega(i)
\end{equation}
The global metric depends on the set of pairs of expected losses for all actions $\{l_{0,0},l_{0,1},l_{1,0},l_{1,1},\cdots\,l_{|\mathcal{A}|,0},l_{|\mathcal{A}|,1}\}$, where $|\mathcal{A}|$ is the cardinality of set $\mathcal{A}$, or, equivalently, on the concave function $l^*$. \figref{fig:global_general} shows an example with $|\mathcal{A}|=4$ actions available. The binary case is when only $|\mathcal{A}|=2$ actions are available: doing-nothing, accepting the risk of paying cost $C_F$ if the system is not working, with $A=0$, or repairing the system at cost $C_R$, with $A=1$ (i.e., $p_{0,0,0}=1, p_{0,1,0}=0, C_0=0, C_1=C_R$). As shown in \figref{fig:global_one}, this setting is defined by $l_{0,0}=0, l_{0,1}=C_F, l_{1,0}=l_{1,1}=C_R$, and the corresponding normalized regret function $rg/C_R$ is bi-linear, with peak $(1-\tilde{p})$ at $p_u=\tilde{p}=C_R/C_F$. 

\subsection{Local metric}
\label{sub:local_metric}

The local metric refers to actions at component level, whose effect depends on components' state $\mathit{s}$. For this approach, we define each action $A$ as a vector $\{a_1, a_2, \cdots, a_N\}$ of $N$ binary entries, where $a_i=1$ if the agent replaces $c_i$, and $a_i=0$ otherwise. Hence the cardinality of the action space is $|\mathcal{A}|=2^N$. 
We assume that the components' repairs are perfect so that transition function $p_{\mathit{s'}|\mathit{s},A}$ is defined as follows: in the vector $s'=[s_1',s_2',\cdots, s_N']$ of states after maintenance, $s_i'=1$ if $a_i=1$, and $s_i'=s_i$ if $a_i=0$.
Function $\mathcal{L}_\RNum{1}(\phi(s'))$ is defined as in \secref{sub:general_voi}, while $\mathcal{L}_\RNum{2}(A)=C_R^\top\cdot A$, with replacing cost vector $C_R=[C_{R,1}, C_{R,2}, \cdots, C_{R,N}]^\top$, where $C_{R,i}$ is the cost of replacing $c_i$.
This model assumes that the accumulated cost is the sum of repair costs for individual component. Other cost models, assuming a more complex cost interaction among components' costs can also be implemented.

After the inspection, the agent selects the optimal subset of components to repair. 
When the inspection outcome is $y_i=c$, the corresponding posterior expected loss, is:
\begin{equation}
    L^\text{L}_{\omega|y_i=c}=\mathbb{E}_{s|y_i=c}\min_A \mathbb{E}_{s'|s,A}\mathcal{L}(s',A)
\end{equation}
Following \eqref{eq:posterior_loss}, the corresponding expected posterior loss is:
\begin{equation}
    \label{eq:local_i}
    L_\omega^\text{L}(i) = (1-h_i)L_{\omega|\mathit{y}_i=1}^\text{L} + h_iL_{\omega|\mathit{y}_i=0}^\text{L}
\end{equation}
and the VoI according to the local metric is $\text{VoI}_\text{L}(i) = L_\pi^\text{L} - L_\omega^\text{L}(i)$, where prior loss $L_\pi^\text{L}$ is computed as in \eqref{eq:prior_loss}.

\subsection{Connection between local and global metrics}
\label{sub:connections_metric}

The local and global metrics refer to different problem classes, which are not nested one into another (given the restrictive rules we impose to the local metric). As the two metrics lead to different analytical approaches, it is natural to ask for a clarifying connection between them. For both metrics, following \secref{sub:general_voi}, one can define a concave function $l^*$ on the belief $p_\mathit{s}$ of the joint condition state $\mathit{s}$ of all components, i.e. on a $N$ dimensional domain (with the linear constraint that $\sum_{i}p_\mathit{s}(i)=1$). Only under the assumption of the global metric, this function can be transformed into a univariate function of the system failure probability $p_u$, as illustrated in \secref{sub:global_metric}. we will not mention that function when dealing with the local metric, while it will be a useful concept for analyzing the properties of the global metric.

\section{Metric properties and inspection priorities on typical networks}
\label{sec:optimal_inspection}

\subsection{Nested posterior intervals for global metric}
\label{sub:priority_for_nested}
\begin{figure}[ht]
    \centering
    \includegraphics[width=.7\linewidth]{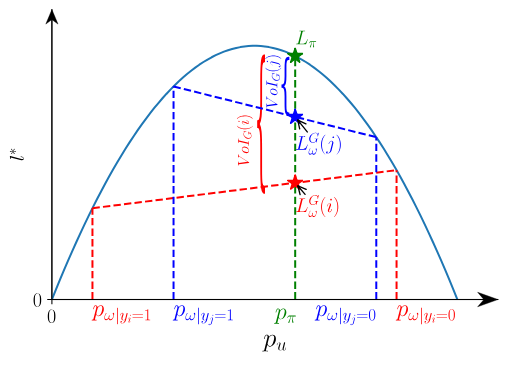}
    \caption{Example of expected loss for the global metric, with nested posterior intervals.}
    \label{fig:Global metric with nested posterior probabilities}
\end{figure}

As we discussed in \secref{sub:global_metric}, the global metric adopts a univariate concave function $l^*$, or $rg$, of $p_\mathit{u}$. An example of such function is shown in \figref{fig:Global metric with nested posterior probabilities}, which can also be interpreted as regret, since it is zero at the limits of the probability domain. Inspecting every component $c_i$, the posterior system failure probability after an alarm is higher than the prior one, and that is in turn higher than the posterior one after a silence: $p_{\omega|y_i=1}\leq p_\pi\leq p_{\omega|y_i=0}$.

Now consider two components $c_i$ and $c_j$. Suppose that a silence on $c_i$ is more reassuring than one on $c_j$ and an alarm from $c_i$ is more worrying than one from $c_j$, i.e. $p_{\omega|\mathit{y}_i=1}\leq p_{\omega|\mathit{y}_j=1}$ and $p_{\omega|\mathit{y}_i=0}\geq p_{\omega|\mathit{y}_j=0}$. Then, for any concave function $l^*$ (or $rg$), the posterior loss inspecting $c_i$ is lower than that inspecting $c_j$ and the VoI of inspecting $c_i$ is higher than that of inspecting $c_j$, i.e. $L^\text{G}_\omega(i)\leq L^\text{G}_\omega(j)$ and $\text{VoI}_\text{G}(i)\geq \text{VoI}_\text{G}(j)$. The proof of this implication is intuitive, examining \figref{fig:Global metric with nested posterior probabilities}, and it is formally given in \appref{app:nested}.

We can also reformulate the implication in terms of “posterior intervals”. Let us define the posterior interval for $c_i$ as $I_i=[p_{\omega|\mathit{y}_j=1}, p_{\omega|\mathit{y}_j=0}]$. If that contains the corresponding interval for $c_j$, i.e. if $I_i\supseteq I_j$, then $\text{VoI}_\text{G}(i)\geq \text{VoI}_\text{G}(j)$. Hence, the importance ranking is invariant respect to the choice of $l^*$, and all possible global metrics favorite one component to inspect respect to another,
which is consistent with Birnbaum's measure \citep{birnbaum1968importance} defined in \secref{sec:introduction}, i.e. $I_i\supseteq I_j\Rightarrow \text{BM}(i)\geq \text{BM}(j)$. However, the reverse implication is not guaranteed and Birnbaum's measure is not necessarily consistent with the global metric.

Moreover, if the posterior intervals are not nested, one can always find a pair of loss functions $\{l^*_\alpha,l^*_\beta\}$, so that $c_i$ has a higher VoI than $c_j$ under ${l^*_\alpha}$, but a lower VoI under ${l^*_\beta}$. To prove that, it suffices to refer to the bi-linear loss function plotted in \figref{fig:expected_loss_rg}. If probability $\tilde{p}$ is not in posterior interval $I_i$ (i.e., $I_i$ is on one side of $\tilde{p}$), then the corresponding VoI, $\text{VoI}_\text{G}(i)$, is zero, as the loss function is linear in that range. If intervals $I_i$ and $I_j$ are not nested, we can find two disjoint intervals: interval $I_{i\diagdown j}$ belongs to $I_i$ but not to $I_j$, interval $I_{j\diagdown i}$ belongs to $I_j$ but not to $I_i$. If $\tilde{p}$ if in $I_{i\diagdown j}$, then $\text{VoI}_\text{G}(i)\geq \text{VoI}_\text{G}(j)=0$, while if $\tilde{p}$ if in $I_{j\diagdown i}$, then $\text{VoI}_\text{G}(j)\geq \text{VoI}_\text{G}(i)=0$. This shows that, for not nested posterior intervals, the priority depends on the adopted loss function.

\subsection{Global metric for parallel systems}
A parallel system will function if at least one of its components is intact. For such systems, the global metric will always give the highest priority to the most reliable component, i.e. the component with the lowest marginal failure probability, independent of the specific loss function $l^*$ adopted, when the inspection quality is the same for all components.
The proof is simple for the special case of perfect sensors, when $\epsilon_{\text{FA}}$ and $\epsilon_{\text{FS}}$ are zero. In that case, if a silence is detected on any component, then posterior system failure probability is zero. As the failure of the system implies the failure of all components, after an alarm on component $c_i$, $p_u$ becomes $p_{\omega|\mathit{s}_j=0}=p_\pi/p_i$. Hence, if $p_i\leq p_j$, then $I_i\supseteq I_j$ and, according to the rule illustrated in \secref{sub:priority_for_nested}, we conclude that $\text{VoI}_\text{G}(i)\geq \text{VoI}_\text{G}(j)$.

When sensors are imperfect, the proof still holds based on Bayes' formula (i.e., on the ratio between joint and marginal probabilities). After a silence on $c_i$, $p_u$ becomes:
\begin{equation}
    \label{eq:silParGlo}
    p_{\omega|\mathit{y}_i=1} = \frac{p_\pi\epsilon_{\text{FS}}}{1-h_i} = \frac{p_\pi\epsilon_{\text{FS}}}{1-\epsilon_{\text{FA}} - K p_i}
\end{equation}
where the second identity follows from \eqref{eq:probAlarm}, and we note again that $K$ is strictly positive. The corresponding probability after an alarm is:
\begin{equation}
    \label{eq:alaParGlo}
    p_{\omega|\mathit{y}_i=0} = \frac{p_\pi(1-\epsilon_{\text{FS}})}{h_i} = \frac{p_\pi(1-\epsilon_{\text{FS}})}{\epsilon_{\text{FA}} + K p_i}
\end{equation}
The denominator of \eqref{eq:silParGlo} is monotonically decreasing with $p_i$, and that of \eqref{eq:alaParGlo} is monotonically increasing with it. Hence, as in the case of perfect sensors, if $p_i\leq p_j$, then $I_i\supseteq I_j$ and, thus, $\text{VoI}_\text{G}(i)\geq \text{VoI}_\text{G}(j)$.

In summary, the ranking of importance measures follows the opposite of the marginal failure probability of the components (i.e., it follows the components' reliability). Hence, in a parallel system, the component, $c_{i^*}$, with highest VoI is the most reliable one. Such result holds for any interdependence among components' states, that is for any distribution $p_\mathit{s}$, when the inspection quality, defined by parameters $\epsilon_{\text{FA}}$ and $\epsilon_{\text{FS}}$, is the same for all components.

\subsection{Global metric for series systems}
A series system works only if all components function properly. In that case, the global metric always prioritizes the most vulnerable component, i.e. that with the highest prior failure probability, regardless of the adopted function $l^*$ or the interdependence among components. The proof is similar to that related to parallel systems. Let us start with the case of perfect sensors. The posterior system failure probability will become $1$ after an alarm on any component, and will become $p_{\omega|\mathit{s}_{j=1}}=1-(1-p_\pi)/(1-p_i)$ after a silence on component $c_i$, which is monotonically increasing with marginal component failure probability $p_i$. Hence the most vulnerable component should be inspected.

For imperfect sensors, after a silence on $c_i$, $p_u$ is (again using \eqref{eq:probAlarm}):
\begin{equation}
    \label{eq:silSerGlo}
    p_{\omega|\mathit{y}_i=1} = 1-\frac{(1-p_\pi)(1-\epsilon_{\text{FA}})}{1-h_i} = 1-\frac{(1-p_\pi)(1-\epsilon_{\text{FA}})}{1-\epsilon_{\text{FA}}-K p_i}
\end{equation}
After an alarm, that probability is:
\begin{equation}
    \label{eq:alaSerGlo}
    p_{\omega|\mathit{y}_i=0} = 1-\frac{(1-p_\pi)\epsilon_{\text{FA}}}{h_i} = 1-\frac{(1-p_\pi)\epsilon_{\text{FA}}}{\epsilon_{\text{FA}} + K p_i}
\end{equation}
The denominator of the fraction in \eqref{eq:silSerGlo} is monotonically decreasing with $p_i$, and that of \eqref{eq:alaSerGlo} is monotonically increasing with it. Hence, if $p_i\geq p_j$, then $I_i\supseteq I_j$ and, thus, $\text{VoI}_\text{G}(i)\geq \text{VoI}_\text{G}(j)$, as in the case of perfect sensors. So, in a series system, regardless the interdependence among components, the inspection ranking follows the marginal component failure probability, and $c_{i^*}$ is the most vulnerable component.

In other words, the most vulnerable component, $c_{i^*}$, is the one to inspect because detecting a silence on that component (i.e. $y_{i^*}=1$) induces the highest reduction of $p_u$, and an alarm (i.e. $y_{i^*}=0$) induces the highest increment in that probability. While the former property is almost trivial, the latter may be less intuitive. After all, $c_{i^*}$ was (relatively) likely to be damaged: so why does an alarm on that component produce the more “surprising” result on the system reliability (compared with alarms on less vulnerable components)? For imperfect inspections, two factors affect the posterior probability. On one hand, after detecting an alarm on $c_{i^*}$, the system can still count on the other components, which are more reliable than $c_{i^*}$ (instead, after an alarm on a safer component, the system can only count on more vulnerable components). Hence, this factor suggests that an alarm of $c_{i^*}$ is less worrying that an alarm on others. However, on the other hand, following Bayes’ rule, an alarm on $c_{i^*}$ produces a relatively high posterior failure probability (at component and at system level), because of the high prior probability that $c_{i^*}$ is damaged. For safer components, the impact of the alarm is diluted by the more optimistic prior information, and the posterior failure probability after an alarm is lower at component level (obviously) and also at system level, as formally proved by \eqref{eq:alaSerGlo}. Hence, this latter factor dominates the former one, and $c_{i^*}$ has the highest VoI. This result depends on the assumption that the sensor accuracy is uniform among components. If the accuracy was higher for a specific component, that component could have the highest VoI, even if it is not the most vulnerable one.

\subsection{Global metric for general systems}
\label{sec:general_system_app}
If the posterior probability interval related to one component is nested with respect to that of another component, then the rule of \secref{sub:priority_for_nested} simply distinguishes the optimal component to inspect. For general systems, the global metric does not always select the most vulnerable or the most reliable component, because the posterior intervals may not be nested, and thus, the rule does not apply.

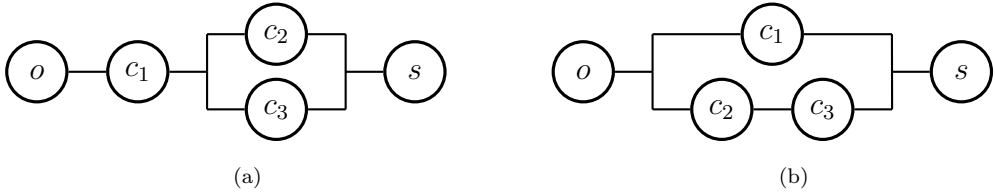
\begin{figure}
    \centering
    \begin{subfigure}[b]{0.47\textwidth}
        \begin{tikzpicture}[
            roundnode/.style={circle, draw=black!100, very thick, minimum width=8mm, inner sep=0pt},
            diamondnode/.style={diamond, draw=black!100, very thick, minimum width=5mm},
            squarenode/.style={rectangle, draw=black!100, very thick, minimum size=7mm},
            ]
            %Nodes
            \node[roundnode] (o) { $o$};
            \node[roundnode, right=0.5cm of o] (n1) { $c_1$};
            \node[coordinate, right=0.5cm of n1] (n1right) {};
            \node[coordinate, above=0.5cm of n1right] (n2left) {};
            \node[roundnode, right=0.5cm of n2left] (n2) { $c_2$};
            \node[coordinate, below=0.5cm of n1right] (n3left) {};
            \node[roundnode, right=0.5cm of n3left] (n3) { $c_3$};
            \node[coordinate, right=0.5cm of n2] (n2right) {};
            \node[coordinate, right=0.5cm of n3] (n3right) {};
            \node[coordinate, below=0.5cm of n2right] (sleft) {};
            \node[roundnode, right=0.5cm of sleft] (s) { $s$};
            
            \draw[thick] (o.east) -- (n1.west);
            \draw[thick] (n1.east) -- (n1right);
            \draw[thick] (n1right) -- (n2left);
            \draw[thick] (n1right) -- (n3left);
            \draw[thick] (n2left) -- (n2.west);
            \draw[thick] (n3left) -- (n3.west);
            \draw[thick] (n2.east) -- (n2right);
            \draw[thick] (n3.east) -- (n3right);
            \draw[thick] (n2right) -- (n3right);
            \draw[thick] (sleft) -- (s.west);
        \end{tikzpicture}
        \caption{}
        \label{fig:3_comps_a}
    \end{subfigure}
    \hfill %
    \begin{subfigure}[b]{0.47\textwidth}
        \begin{tikzpicture}[
            roundnode/.style={circle, draw=black!100, very thick, minimum width=8mm, inner sep=0pt},
            diamondnode/.style={diamond, draw=black!100, very thick, minimum width=5mm},
            squarenode/.style={rectangle, draw=black!100, very thick, minimum size=7mm},
            ]
            %Nodes
            \node[roundnode] (o) { $o$};
            \node[coordinate, right=0.5cm of o] (oright) {};
            \node[coordinate, above=0.5cm of oright] (n1left) {};
            \node[coordinate, below=0.5cm of oright] (n2left) {};
            \node[roundnode, right=0.5cm of n2left] (n2) { $c_2$};
            \node[roundnode, right=0.5cm of n2] (n3) { $c_3$};
            \node[roundnode, right=1.17cm of n1left] (n1) { $c_1$};
            \node[coordinate, right=1.17cm of n1] (n1right) {};
            \node[coordinate, right=0.5cm of n3] (n3right) {};
            \node[coordinate, above=0.5cm of n3right] (sleft) {};
            \node[roundnode, right=0.5cm of sleft] (s) { $s$};
            
            \draw[thick] (o.east) -- (oright);
            \draw[thick] (oright) -- (n1left);
            \draw[thick] (oright) -- (n2left);
            \draw[thick] (n1left) -- (n1.west);
            \draw[thick] (n1.east) -- (n1right);
            \draw[thick] (n2left) -- (n2.west);
            \draw[thick] (n2.east) -- (n3.west);
            \draw[thick] (n3.east) -- (n3right);
            \draw[thick] (n1right) -- (n3right);
            \draw[thick] (sleft) -- (s.west);
        \end{tikzpicture} 
        \caption{}
        \label{fig:3_comps_b}
    \end{subfigure}
    \caption{Series-parallel (a) and parallel-series (b) 3-component system.}
    \label{fig:3_comps}
\end{figure}

We illustrate this discussing two simple examples of 3-component systems, with perfect sensors. \figref{fig:3_comps_a} shows a system where component $c_1$ is in series with the parallel subsystem made up by components $c_2$ and $c_3$.  Intuitively, component $c_1$ should be inspected, as it is a “bottleneck” of the system, and so it seems topologically more important. Detecting that $c_1$ is not working takes $p_u$ to one, i.e. to a value higher than that related to an alarm on $c_2$ or on $c_3$. A silence detected on $c_1$ takes $p_u$ to the joint failure probability, $p_{\omega|s_1=1}$, which is determined by components $c_2$ and $c_3$, and is $(p_2p_3)$ if they are independent. Instead, a silence detected on $c_2$ (or on $c_3$) takes $p_u$ to $p_{\omega|s_2=1}=p_1$. Hence, posterior interval $I_1$ contains the other two if $p_{\omega|s_1=1}$ is less than $p_1$. On the contrary, if $p_1$ is less than $p_{\omega|s_1=1}$, the posterior intervals are not nested, and the priority depends on the selected loss function $l^*$. This result confirms the intuition that if $c_1$ it much safer than the other components, it may not be the most important component to inspect (trivially, if $p_1$ is zero while $p_2$ and $p_3$ are positive, then $c_1$ has the lowest priority).

In the example of \figref{fig:3_comps_b}, component $c_1$ is parallel with a series subsystem made up by components $c_2$ and $c_3$. Again, $c_1$ seems topologically more important. After a silence on $c_1$, $p_u$ is zero, a value lower than that related to silence on $c_2$ or on $c_3$. An alarm on $c_1$ takes $p_u$ to $1-r_{2,3}$, where $r_{2,3}$ is the joint survival probability of the other two components, that is $(1-p_2)(1-p_3)$ for independent components, while an alarm on $c_2$ (or on $c_3$) takes $p_u$ to $p_1$. Hence, posterior interval $I_1$ contains the others if $p_1$ is less than $1-r_{2,3}$ i.e., for independent components, if $p_1$ is less than $p_2+p_3-p_2p_3$. Approximating this latter value with $p_2+p_3$, we conclude that the global metric gives higher priority to $c_1$ when $p_1$ is less than $(p_2+p_3)$. If $p_1$ is higher than that, priority depends on the selected loss function $l^*$. This confirms the intuition that, if $c_1$ is much more vulnerable than the other components, it is better to inspect others (in the limit case where $p_1$ is one, $\text{VoI}_\text{G}(1)$ is zero). These two examples illustrate how the topological role of a component matters, but also its failure probability: in some schemes a high failure probability guarantees a high priority, in others a low probability does.

\begin{figure}
    \centering
    \begin{tikzpicture}[
        roundnode/.style={circle, draw=black!100, very thick, minimum width=10mm, inner sep=0pt},
        diamondnode/.style={diamond, draw=black!100, very thick, minimum width=5mm},
        squarenode/.style={rectangle, draw=black!100, very thick, minimum size=7mm},
        ]
        %Nodes
        \node[roundnode] (o) {\Large $o$};
        \node[coordinate, right=0.5cm of o] (b1) {};
        \node[coordinate, above=1cm of b1] (b2) {};
        \node[coordinate, below=1cm of b1] (b3) {};
        \node[roundnode, right=0.5cm of b2] (f1) {\small $0.5\%$};
        \node[coordinate, right=0.5cm of f1] (fb1) {};
        \node[coordinate, above=0.8cm of fb1] (fb2) {};
        \node[coordinate, below=0.8cm of fb1] (fb3) {};
        \node[roundnode, right=0.5cm of fb2] (c11) {\small$2.3\%$};
        \node[roundnode, right=0.5cm of fb3] (c21) {\small$2.3\%$}; 
        \node[roundnode, right=0.5cm of c11] (c1) {\small$1\%$};
        \node[roundnode, right=0.5cm of c21] (c2) {\small$20\%$};
        \node[coordinate, right=0.5cm of c1] (fb4) {};
        \node[coordinate, right=0.5cm of c2] (fb5) {};
        \node[coordinate, below=0.8cm of fb4] (fb6) {};
        \node[coordinate, right=0.5cm of fb6] (fb7) {};

        \node[roundnode, right=2.5cm of b3] (f2) {\small$90\%$};
        \node[coordinate, below=2cm of fb7] (f2b) {};
        \node[coordinate, below=1cm of fb7] (fend) {};
        \node[roundnode, right=0.5cm of fend] (s) {\Large $s$};

        \node [above=0.03cm of c1] { $c_1$};
        \node [above=0.03 cm of c2] { $c_2$}; 

        \draw [thick] (o.east) -- (b1.west);
        \draw [thick] (b1.north) -- (b2.south);
        \draw [thick] (b1.south) -- (b3.north);
        \draw [thick] (b2.east) -- (f1.west);
        \draw [thick] (f1.east) -- (fb1.west);
        \draw [thick] (fb1) -- (fb2);
        \draw [thick] (fb1) -- (fb3);
        \draw [thick] (fb2) -- (c11.west);
        \draw [thick] (fb3) -- (c21.west);
        \draw [thick] (c11.east) -- (c1.west);
        \draw [thick] (c21.east) -- (c2.west);
        \draw [thick] (c1.east) -- (fb4);
        \draw [thick] (c2.east) -- (fb5);
        \draw [thick] (fb4) -- (fb6);
        \draw [thick] (fb5) -- (fb6);
        \draw [thick] (fb6) -- (fb7);
        \draw [thick] (b3) -- (f2.west);
        \draw [thick] (f2.east) -- (f2b);
        \draw [thick] (f2b) -- (fb7);
        \draw [thick] (fend) -- (s.west);            
    \end{tikzpicture}
    \caption{Example of a system where posterior probability intervals for components $c_1$ and $c_2$ are not nested}
    \label{fig:general_nest_counter}
\end{figure}
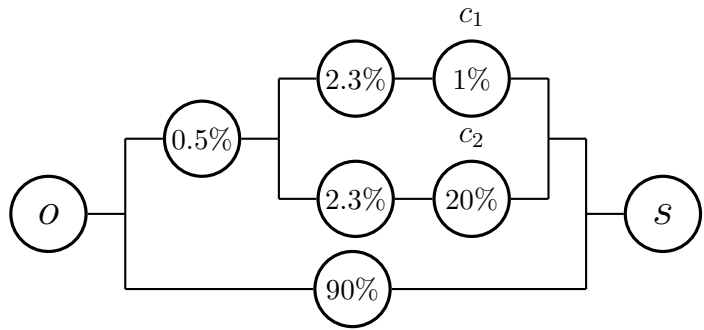

\begin{figure}[h]
    \centering
    \includegraphics[width=.5\textwidth]{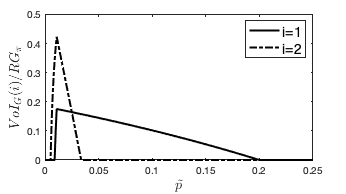}
    \caption{Normalized VoI depending on peak probability $\tilde{p}$.}
    \label{fig:bilinear}
\end{figure}

We discuss now a more general example, focusing on two components, $c_1$ and $c_2$. The components’ role is completely described by the system failure probability for each of the $2^2$ joint conditions of the pair of components, that we assume as $p_{\omega|s_1=1,s_2=1}=0.5\%$, $p_{\omega|s_1=1,s_2=0}=p_{\omega|s_1=0,s_2=1}=2.5\%$, $p_{\omega|s_1=0,s_2=0}=90\%$ (so the roles played by the two components are the same). We also assume that $p_1=1\%$ and $p_2=20\%$ (so that the $c_2$ is significantly more reliable than $c_1$), the states of all components are independent, and inspections are perfect (i.e., $y_i$=$s_i$). \figref{fig:general_nest_counter} shows the diagram of a system consistent with these values. The interval of posterior probabilities $I_i$ is $[0.90\%,20.0\%]$ for $i=1$ and $[0.52\%,3.38\%]$ for $i=2$, while $p_\pi$ is $1.09\%$ (these results are directly related to the assumed values, e.g. $p_{\omega|y_2=1}=p_{\omega|s_1=0,s_2=1}p_1+p_{\omega|s_1=1,s_2=1}(1-p_1)$). The intervals are not nested, hence the rule in \secref{sub:priority_for_nested} does not apply, and the VoI depends on the specific function $l^*$. \figref{fig:bilinear} refers to the bi-linear regret function for binary actions plotted in \figref{fig:expected_loss_rg}, and mentioned in Section 2.2, with peak at $\tilde{p}$, and it shows how the VoI related to each component, normalized by prior regret $rg_\pi$, varies as a function of $\tilde{p}$. When $\tilde{p}$ is below $p_{\omega|y_2=1}=0.52\%$ (i.e., when $C_R$ is below $0.52\%$ of $C_F$), the VoI of each component is nil, as the posterior decision is always to repair. Then, $\text{VoI}_G(2)$ increases up to about $42\%$ of $RG_\pi$ when $\tilde{p}=p_\pi$ (i.e., for that condition observing $y_2$ is worth $42\%$ of the value of observing $u$), then it decreases down to zero at $\tilde{p}=p_{\omega|y_2=0}=3.38\%$. For $\tilde{p}$ higher than that value, the posterior decision is always to accept the risk. The behavior of $\text{VoI}_\text{G}(1)$ is similar: it is zero outside $I_1$, and it peaks at $p_\pi$, where it is about $17\%$ of $RG_\pi$. Clearly, the optimal inspection decision depends on $\tilde{p}$, i.e. on the decision-making problem shaping function $l^*$: as apparent in \figref{fig:bilinear}, if the repair cost is cheaper than $2.5\%$ of the cost of failure, it is more convenient to inspect the more reliable $c_2$, while it is better to inspect the less reliable $c_1$ for an higher repair cost.

% \begin{figure}[h]
%     \centering
%     \includegraphics[width=.7\textwidth]{figure/example_general_nest.png}
%     \caption{General system example with posterior probability intervals ($I_1$ and $I_2$) not nested}
%     \label{fig:general_nest_counter}
% \end{figure}

\subsection{Local metric on parallel systems}
The local metric, as defined in \secref{sub:local_metric}, will select the most reliable component in a parallel system, consistently with the global metric. This is because, in a parallel systems, replacing one component guarantees the functioning of the system. Hence, the agent faces a binary decision: doing nothing or repairing the less expensive component, at cost $\min_iC_{R,i}$, and this problem setting is equivalent to that of the global metric, with the bi-linear function $l^*$ of \figref{fig:global_one}.
Therefore, the local and global metric have identical conclusions about the optimal inspection. 

\subsection{Local metric on series systems}
\begin{figure}
    \centering
    \includegraphics[width=.7\textwidth]{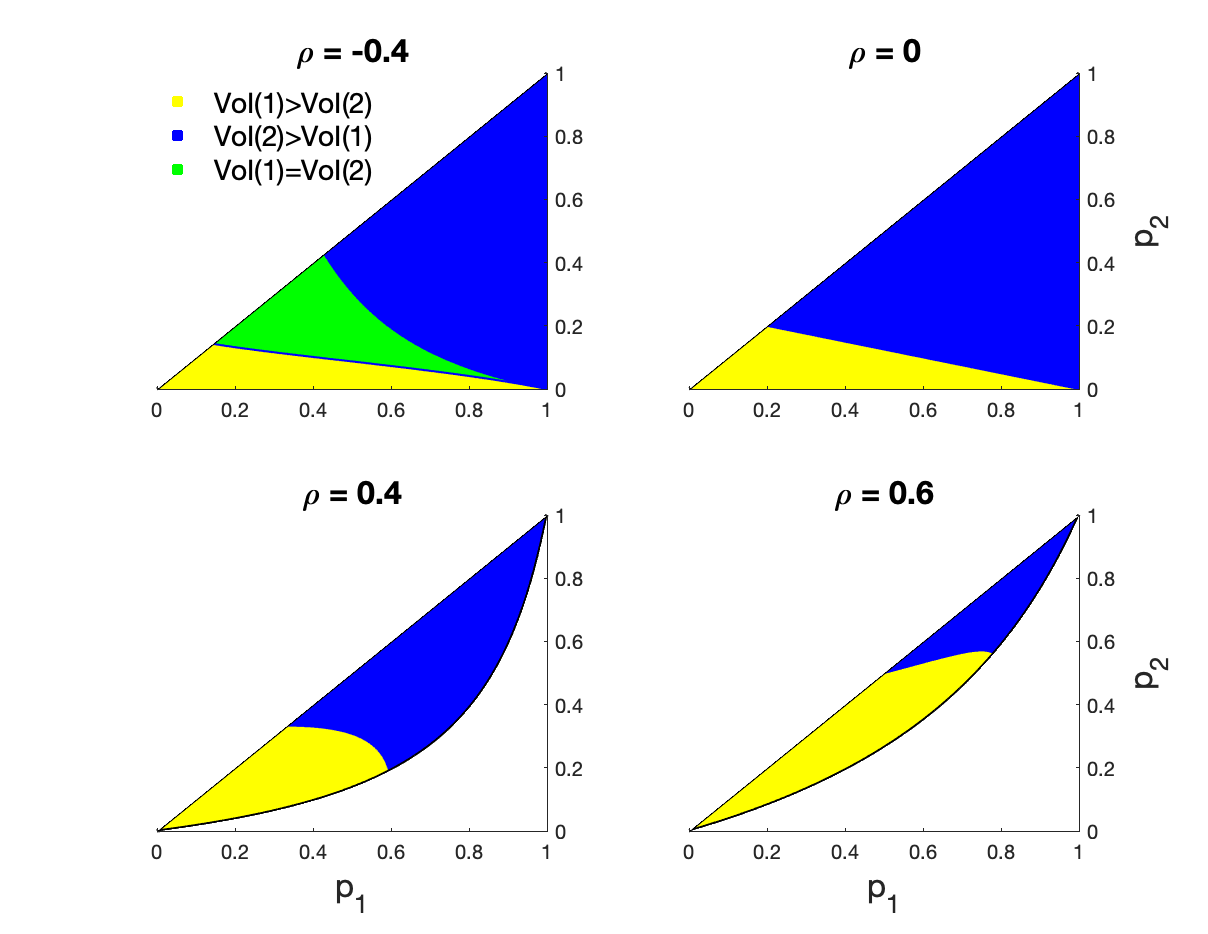}
    \caption{Optimal inspection $i^*$ when two components are dependent, depending on correlation $\rho$.}
    \label{fig:Optimal inspections when two components are dependent}
\end{figure}

Under the local metric, the optimal component to inspect in a series system is not always the the most vulnerable one, i.e. that identified by the global metric.
We start discussing the case of a system with two components, $c_1$ and $c_2$, with identical repair costs, $C_{R_1}=C_{R_2}=C_R$, and equipped with perfect sensors. Let us also assume that $C_R\leq C_F/2$, so that the cost for replacing both components is less than the failure cost. Hence, if any component $c_i$ is detected as damaged, it is necessary to repair it ($A_i=1$), to avoid paying the failure cost. After the replacement, the system failure probability is the posterior failure probability of the uninspected component, and that should be also replaced if the corresponding risk is above the repair cost, so that the posterior expected maintenance cost for that component is $R(i,x)=\min\{C_R,p_{\omega|s_i=x,A_i=1-x}C_F\}$, with $x=0$. Instead, if the inspected component works, it has not to be repaired, and the state of the uninspected one is decided by comparing repair cost and system failure risk, so that the expected posterior cost is $R(i,1)$. Hence, the expected posterior loss is $L^\text{L}_{\omega}(i)=p_iC_R+p_iR(i,0)+(1-p_i)R(i,1)$. In the special case of independent components, for any outcome $x$, probability $p_{\omega|s_i=x}$ is identical to the prior failure probability, $p_j$, of the uninspected component $c_j$, so that $R(i,0)=R(i,1)$, and $L^\text{L}_{\omega}(i)=p_iC_R+\min\{C_R,p_jC_F\}$. If we refer to bi-linear regret function $rg$ of \figref{fig:3_comps_b}, we conclude that, for each component $c_i$, $\text{VoI}_\text{L}(i)=\text{RG}_\omega(i)$, and so we should inspect the component with the higher value of $\text{RG}_\omega(i)$. If both prior failure probabilities are below $\tilde{p}=C_R/C_F$, the local metric will prioritize the more vulnerable component. However, the failure probability of a component is above $\tilde{p}$, the higher that probability is, the lower the corresponding VoI. \figref{fig:Optimal inspections when two components are dependent} shows the optimal inspection policies for $\tilde{p}=0.2$, $p_1\geq p_2$, and different correlation coefficient $\rho$ between variables $s_1$ and $s_2$. We have discussed the case when $\rho$ is zero. When it is positive, the domain of feasible pairs ($p_1$,$p_2$) shrinks but, inside the feasible domain, the region where it is more convenient to inspect the more vulnerable component expands. When the correlation is negative, for any feasible pair $\{p_1,p_2\}$, the VoI is the same for both components.

When the series system is made up by a higher number of components, we can provide a simple approximation if their states are independent and the failure probabilities are relatively low. In that case, the risk $\mathbb{E}[\mathcal{L}_I]$ can be approximated as that of a “cumulative system”\citep{Malings2016}. In a cumulative system, individual costs are associated to the failure of each component, and accumulated to get the system-level cost (hence, no component-to-system function $\phi$ is defined for these systems).
To show this, we recall that, for a series system with independent components, the risk is:
\begin{equation}
    \mathbb{E}_\text{ser.}[\mathcal{L}_I]=C_F[1-\prod_i(1-p_i)^{1-A_i}]
\end{equation}
For a cumulative system with component failure cost $C_F$, it is:
\begin{equation}
   \mathbb{E}_\text{cum.}[\mathcal{L}_I]=C_F\sum_i(1-A_i)p_i
\end{equation}
By linearizing the former expression (so neglecting higher order terms), the two risks becomes identical. For a cumulative system, it is straightforward to evaluate the benefit of inspecting component $c_i$ as that related to the selection of action $A_i$ and, when sensors are perfect, it is $\text{VoI}_\text{L}(i)=\text{RG}_\omega(i)$ (for imperfect sensors, one has to subtract the posterior regret).
These results are also consistent with the case discussed above, where $N=2$.

\section{Computational complexity and Heuristic}
\label{sec:heuristic}

\subsection{Complexity of VoI computation}
The computational difficulty of solving \eqref{eq:ystar} varies with different metrics, but is generally intimidating for large systems. The core step of the computational process is solving the reliability problem, identifying the system failure probability $p_u$, depending on actions and observations. This analysis is nested into the optimization of the maintenance actions.

To compute the risk $\mathbb{E}[\mathcal{L}_\RNum{1}]$, one has to assess the network connectivity for each of the $2^N$ components' states.
\cite{Song2003BoundsReliability} proposed matrix-based method to compute system reliability based on a components' condition matrix with each row representing one state, and a binary condition vector with each entry representing whether the system is functioning or not at that specific components' state. The general computation complexity is $\mathcal{O}(N \times 2^N)$. 
When the joint distribution of the components' states is known, the computational complexity of system reliability is linear.
An approximate estimation can also be achieved based on Monte Carlo simulations \cite{kroese2013handbook}.

For the global metric, $\mathbb{E}[\mathcal{L}_\RNum{2}]$ can be determined in $\mathcal{O}(1)$ time once we have computed the posterior system failure probabilities. Therefore, to select the component with highest VoI among $N$ components will cost $\mathcal{O}(N\times 2^N)$.

For the local metric, there is an additional computation step before assessing the VoI. $\mathbb{E}[\mathcal{L}_\RNum{2}]$ is optimized among $2^N$ combinations of maintenance actions. Suppose that, based on different inspection outcomes, the agent can select an arbitrary subset of the components to replace/repair, the computation is generally $\mathcal{O}(N\times 2^N\times 2^N) =\mathcal{O}(N\times 2^{2N} )$.

\subsection{Approximation for local metric}
In this section, we propose a simple heuristic approach for approximating the local metric, to reduce the computational complexity related to the optimization of maintenance actions depending on the inspection outcome. 

Let us define $A_\pi=\{a_{\pi,1},a_{\pi,2},\cdots,a_{\pi,N}\}$ as the prior maintenance plan, $A_\omega=\{a_{\omega,1}, a_{\omega,2},\cdots,a_{\omega,N} \}$ as the posterior one, and $L_\pi$ is the prior optimal loss related to $A_\pi$, as defined in \secref{sub:general_voi}. We assume that $A_\pi$ and $L_\pi$ have been identified. Consider inspecting component $c_i$. The proposed heuristic assumes that the agent confirms all actions for uninspected components (i.e., $\forall j\neq i, a_{\omega, j}=a_{\pi,j}$). Only the posterior action on the inspected component, $a_{\omega, i}$, depends on the inspection’s outcome, $y_i$. If the prior action for $c_i$ is to do-nothing (i.e., if $a_{\pi,i}=0$) and the inspection’s outcome is silence (i.e., if $y_i=1$), or if the prior action is to repair (i.e., if $a_{\pi,i}=1$) and the inspection produces an alarm (i.e., if $y_i=0$), then the agent will confirm the prior action also for the inspected component (i.e., if $y_i\neq a_{\pi,i}$, then $A_\omega=A_\pi$). Instead, if an alarm is detected on a prior unrepaired component, or if a silence is detected on a prior repaired component (i.e., if $y_i=a_{\pi,i}$), then the agent consider the two alternatives: to repair $c_i$ or not to. One of the two alternatives is, again, to completely confirm the prior plan (i.e. $A_\omega=A_\pi$), and thus the prior loss $L_\pi$ associated with this option is already known. The agent computes the expected cost of the alternative plan (where only action $a_{\omega,i}$ is reversed), and executes the best option, i.e. that related to the minimum expected cost. The computational saving is related to the avoidance of the full posterior optimization in set $\mathcal{A}$.

One argument supporting the choice of this heuristic is that it is consistent with the optimal behavior in some special cases. For example, when the high risk forces the agent to be conservative. To model that, suppose that (i) the prior decision is to do-nothing (i.e., $\forall i, a_{\pi, i}=0$), that (ii) a detected silence cannot increase the system failure probability (i.e. $\forall i, p_{\omega|y_i=1}\leq p_\pi$), that (iii) the do-nothing option is still optimal when the probability of failure decreases and that (iv) a component sending an alarm must be replaced, as its posterior failure probability is too high to be tolerated. Condition (iii) is not obviously satisfied even if the first two are, as the prior decision might also be doing-nothing for another reason, i.e. because the agent is pessimistic about the components’ condition. For such a pessimistic agent, it is not worth to repair any set of components: repairing few components may be ineffective, and repairing many components may be too expensive. However, detecting a functioning component may suggest the pessimistic agent to invest in the system’s repairing, by replacing other components. Condition (iii) forbids the occurrence of this process, by assuming agent’s optimism about the system condition. To prove that the heuristic is optimal under conditions (i-iv), we must show that the optimal response to an alarm on component $c_i$ cannot be to repair any other component. Because of (iv), $c_i$ must be repaired. Now suppose that, component $c_j$ is also to be repaired. This implies the following inequality:
\begin{equation}
    C_{R,i}+C_{R,j}+C_F p_{\omega|y_i=0,a_i=1,a_j=1}\leq C_{R,i}+C_Fp_{\omega|y_i=0,a_i=1}
\end{equation}
If, as assumed before, repairs are perfect and components’ states are independent, then $p_{\omega|y_i=0,a_i=1}=p_{\omega|s_i=1}=p_{\omega|y_i=1}$, and $p_{\omega|y_i=0,a_i=1,a_j=1}=p_{\omega|s_i=1,s_j=1}=p_{\omega|y_i=1,a_j=1}$, so that previous inequality can be re-written, subtracting $C_{R,i}$ to both terms, as:
\begin{equation}
    C_{R,j}+C_F p_{\omega|y_i=1,a_j=1}\leq C_F p_{\omega|y_i=1}
\end{equation}
Indicating that repairing $c_j$ should be the optimal response to a silence on $c_i$, but this violates conditions (i-iii), showing that only $c_i$ should be repaired after an alarm on that component. Of course, if conditions (i-iv) are not satisfied, there is no guarantee that the heuristic is truly optimal.

The VoI defined by the heuristic is certainly non-negative, as the prior maintenance plan can be confirmed, if the collected observations do not suggest any improvement. Moreover, given that the heuristic limits the domain of the posterior actions, the corresponding VoI cannot be higher than the original one assessed by the local metric.

\section{Examples of Network Analysis}
\label{sec:examples}

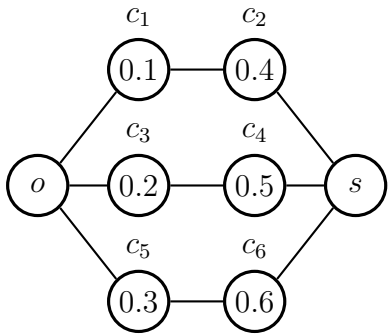
\begin{figure}
    \centering
    \begin{tikzpicture}[
        roundnode/.style={circle, draw=black!100, very thick, minimum width=8mm, inner sep=0pt},
        diamondnode/.style={diamond, draw=black!100, very thick, minimum width=5mm},
        squarenode/.style={rectangle, draw=black!100, very thick, minimum size=7mm},
        ]
        %Nodes
        \node[roundnode] (o) { $o$};
        \node[roundnode, right=0.5cm of o] (n3) {$0.2$};
        \node[roundnode, above=0.7cm of n3] (n1) {$0.1$};
        \node[roundnode, below=0.7cm of n3] (n5) {$0.3$};
        \node[roundnode, right=0.7cm of n1] (n2) {$0.4$};
        \node[roundnode, right=0.7cm of n3] (n4) {$0.5$};
        \node[roundnode, right=0.7cm of n5] (n6) {$0.6$};
        \node[roundnode, right=0.5cm of n4] (s) { $s$};
        \node[above=0.01cm of n1] { $c_1$};
        \node[above=0.01cm of n2] { $c_2$};
        \node[above=0.01cm of n3] { $c_3$};
        \node[above=0.01cm of n4] { $c_4$};
        \node[above=0.01cm of n5] { $c_5$};
        \node[above=0.01cm of n6] { $c_6$};

        \draw[thick] (o.north east) -- (n1.south west);
        \draw[thick] (o.east) -- (n3.west);
        \draw[thick] (o.south east) -- (n5.north west);
        \draw[thick] (n1.east) -- (n2.west);
        \draw[thick] (n3.east) -- (n4.west);
        \draw[thick] (n5.east) -- (n6.west);
        \draw[thick] (n2.south east) -- (s.north west);
        \draw[thick] (n4.east) -- (s.west);
        \draw[thick] (n6.north east) -- (s.south west);
    \end{tikzpicture}
    \caption{Counter-intuitive example 1}
    \label{fig:counter_1}
\end{figure}

\begin{figure}
    \centering
    \begin{subfigure}[b]{0.47\textwidth}
        \includegraphics[width=\textwidth]{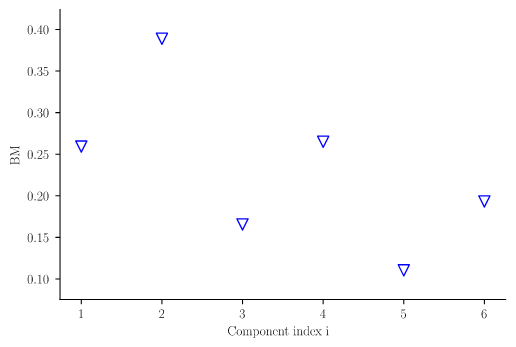}
        \caption{}
        \label{fig:VoI under Birnbaum's measure for network in}
    \end{subfigure}
    \hfill
    \begin{subfigure}[b]{0.47\textwidth}
        \includegraphics[width=\textwidth]{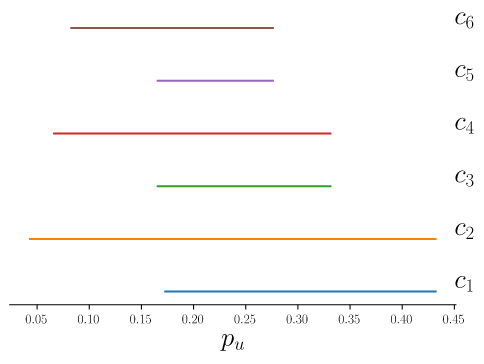}
        \caption{}
        \label{fig:Posterior probabilities for each components}
    \end{subfigure}
    \caption{BM for the network in \figref{fig:counter_1} (a), and corresponding posterior intervals (b).}
    \label{fig:Normalized VoI under global/local/heuristic approaches for network bm}
\end{figure}

We analyze three examples of networks. The first one is the 6-component network in \figref{fig:counter_1}, where the failure probability of each component is listed inside the corresponding node. We start considering perfect inspections and independent components. The corresponding values of the $\text{BM}$s are shown in \figref{fig:VoI under Birnbaum's measure for network in}, and $c_2$ has the highest importance in $\text{BM}$.

\begin{figure}
    \centering
    \begin{subfigure}[b]{0.47\textwidth}
        \includegraphics[width=\textwidth]{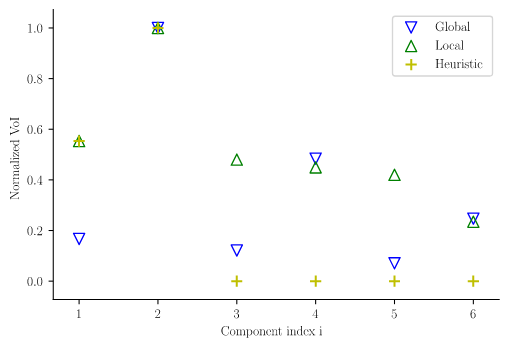}
        \caption{}
        \label{fig:cfcr10}
    \end{subfigure}
    \hfill
    \begin{subfigure}[b]{0.47\textwidth}
        \includegraphics[width=\textwidth]{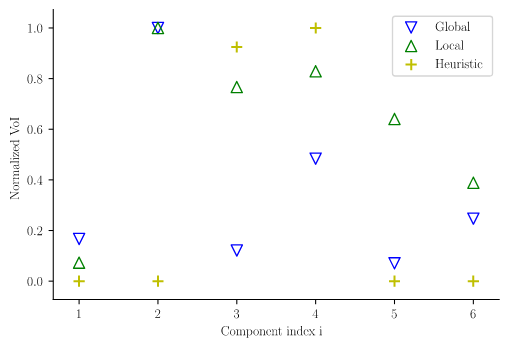}
        \caption{}
        \label{fig:cfcr5}
    \end{subfigure}
    \caption{Normalized VoI for the network in \figref{fig:counter_1}, with $C_F/C_{R_i}=10$ (a), and with $C_F/C_{R_2}=5, C_F/C_{R_i}=10,\forall i\neq 2$ (b).}
    \label{fig:Normalized VoI under global/local/heuristic approaches for network in}
\end{figure}

\figref{fig:Posterior probabilities for each components} shows the posterior probabilities intervals $I_i$ for all components. All intervals are nested in $I_2$, thus component $c_2$ has the highest VoI, according to the global metric, regardless of the loss function we adopt (and so it has also the highest BM, as noted above). We divide the VoI of each component by the maximum VoI among all the components under the same metric to get normalized VoI. The normalized VoI under the global metric, for loss function $l^*(p_u)=p_u(1-p_u)$, is shown in \figref{fig:Normalized VoI under global/local/heuristic approaches for network in}. 

For the local metric, we assume that $C_F/C_{R,i}=10$, for every component $c_i$, i.e. the cost of system failure is ten times the cost of replacing one component. The optimal prior maintenance action is to replace component $c_2$. As shown in \figref{fig:cfcr10}, the local metric and the heuristic both identify $c_2$ as the component with highest VoI. 

However, if the maintenance cost for $c_2$ increases to $C_F/C_{R,2}=5$ while the cost for the others remains the same, the optimal prior action becomes replacing $c_4$. \tabref{tab:Optimal posterior action for the counter-intuitive example 1} reports the optimal posterior actions depending on the inspection outcome, for this new assumption on the costs. As shown in \figref{fig:cfcr5}, the local metric still gives the highest inspection priority to $c_2$ (as the global metric does), but the heuristic selects $c_4$ instead. This is because the posterior optimal action may not include repairing $c_4$ (e.g. after a silence on $c_2$), or it may include repairing uninspected components (e.g. after an alarm on $c_1$, $c_3$ is to be repaired). The heuristic overestimates the value of inspecting $c_4$, and its assessment is not consistent with that of the local metric.

\begin{table}
    \centering
    \begin{tabular}{ccc}
        \toprule
        \diagbox[width=15em]{Inspected component}{Inspection outcome}&
        Silence ($y_i=1$) & Alarm ($y_i=0$)\\
        \midrule
        $c_1$ & $\{c_4\}$ & $\{c_3,c_4\}$\\
        $c_2$ & $\emptyset$ & $\{c_3,c_4\}$ \\
        $c_3$ & $\{c_4\}$ & $\{c_3,c_4\}$ \\
        $c_4$ & $\emptyset$ & $\{c_4\}$ \\
        $c_5$ & $\{c_6\}$ & $\{c_4\}$ \\
        $c_6$ & $\{\emptyset\}$ & $\{c_6\}$ \\
        \bottomrule
    \end{tabular}
    \caption{Posterior subset of components to be repaired for the network in \figref{fig:counter_1}.}
    \label{tab:Optimal posterior action for the counter-intuitive example 1}
\end{table}

\begin{figure}
    \centering
    \begin{subfigure}[b]{0.47\textwidth}
        \includegraphics[width=\textwidth]{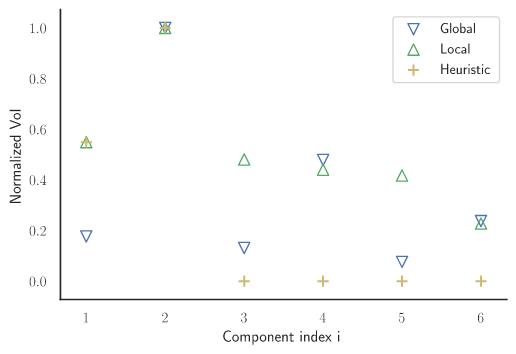}
        \caption{}
        \label{fig:noise_counter_11}
    \end{subfigure}
    \hfill
    \begin{subfigure}[b]{0.47\textwidth}
        \includegraphics[width=\textwidth]{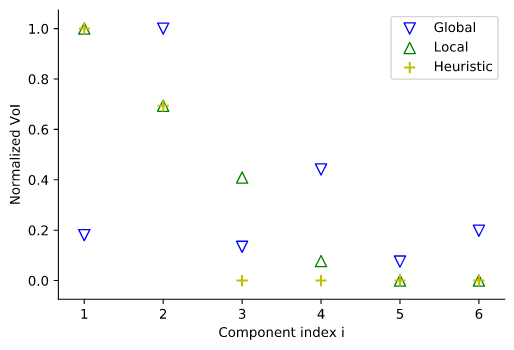}
        \caption{}
        \label{fig:noise_counter_12}
    \end{subfigure}
    \caption{Normalized VoI for the network in \figref{fig:counter_1}, with $\epsilon_{FA}=\epsilon_{FS}=0.01$ (a), and $\epsilon_{FA}=0.01, \epsilon_{FS}=0.40$ (b).}
    \label{fig:VoI for network counter}
\end{figure}

Error rates in imperfect inspections also affect the optimal decision. We now assume, again, that $C_F/C_{R,i}=10$ for every component $c_i$, but now inspections are imperfect: when $\epsilon_{FA}=\epsilon_{FS}=0.01$, the corresponding VoI, shown in \figref{fig:noise_counter_11}, is similar to the perfect inspection case shown in \figref{fig:cfcr10}, and $c_2$ has the highest VoI. But when the type $\RNum{2}$ error rate $\epsilon_{FS}$ is increased to $0.40$, the VoI becomes that shown in \figref{fig:noise_counter_12}, and component $c_1$ gains the highest priority for the local metric and heuristic approach. 

\begin{figure}
    \centering
    \begin{tikzpicture}[
        roundnode/.style={circle, draw=black!100, very thick, minimum width=6mm, inner sep=0pt},
        diamondnode/.style={diamond, draw=black!100, very thick, minimum width=5mm},
        squarenode/.style={rectangle, draw=black!100, very thick, minimum size=7mm},
        ]
        %Nodes
        \node[roundnode] (o) { $o$};
        \node[coordinate, right=0.8cm of o] (oright) {};
        \node[roundnode, above=1cm of oright] (n1) { $c_1$};
        \node[roundnode, below=1cm of oright] (n8) { $c_8$};
        \node[coordinate, right=0.8cm of n1] (n1right) {};
        \node[roundnode, above=0.5cm of n1right] (n2) { $c_2$};
        \node[roundnode, below=0.5cm of n1right] (n3) { $c_3$};
        \node[roundnode, right=0.8cm of n1right] (n4) { $c_4$};
        \node[coordinate, right=0.8cm of n4] (n4right) {};
        \node[roundnode, above=0.5cm of n4right] (n5) { $c_5$};
        \node[roundnode, below=0.5cm of n4right] (n6) { $c_6$};
        \node[roundnode, right=0.8cm of n4right] (n7) { $c_7$};
        \node[coordinate, below=1cm of n7] (sleft) {};
        \node[roundnode, right=0.8cm of sleft] (s) { $s$};
        \node[coordinate, right=0.8cm of n8] (n8right) {};
        \node[roundnode, above=0.5cm of n8right] (n9) { $c_9$};
        \node[roundnode, below=0.5cm of n8right] (n10) { $c_{10}$};
        \node[roundnode, right=0.8cm of n8right] (n12) { $c_{12}$};
        \node[roundnode, above=1cm of n12] (n11) { $c_{11}$};
        \node[roundnode, below=1cm of n12] (n13) { $c_{13}$};
        \node[roundnode, right=1.6cm of n9] (n14) { $c_{14}$};
        \node[roundnode, right=1.6cm of n10] (n15) { $c_{15}$};
        \node[roundnode, right=1.6cm of n12] (n16) { $c_{16}$};
        
        \draw[thick] (o.north east) -- (n1.south west);
        \draw[thick] (o.south east) -- (n8.north west);
        \draw[thick] (n1.north east) -- (n2.south west);
        \draw[thick] (n1.south east) -- (n3.north west);
        \draw[thick] (n2.south east) -- (n4.north west);
        \draw[thick] (n3.north east) -- (n4.south west);
        \draw[thick] (n4.north east) -- (n5.south west);
        \draw[thick] (n5.south east) -- (n7.north west);
        \draw[thick] (n4.south east) -- (n6.north west);
        \draw[thick] (n6.north east) -- (n7.south west);
        \draw[thick] (n7.south east) -- (s.north west);
        \draw[thick] (n8.north east) -- (n9.south west);
        \draw[thick] (n8.south east) -- (n10.north west);
        \draw[thick] (n9.north east) -- (n11.south west);
        \draw[thick] (n9.south east) -- (n12.north west);
        \draw[thick] (n10.north east) -- (n12.south west);
        \draw[thick] (n10.south east) -- (n13.north west);
        \draw[thick] (n11.south east) -- (n14.north west);
        \draw[thick] (n12.north east) -- (n14.south west);
        \draw[thick] (n12.south east) -- (n15.north west);
        \draw[thick] (n13.north east) -- (n15.south west);
        \draw[thick] (n14.south east) -- (n16.north west);
        \draw[thick] (n15.north east) -- (n16.south west);
        \draw[thick] (n16.north east) -- (s.south west); 
    \end{tikzpicture}
    \caption{Example of 16 component network}
    \label{fig:general_network_example}
\end{figure}
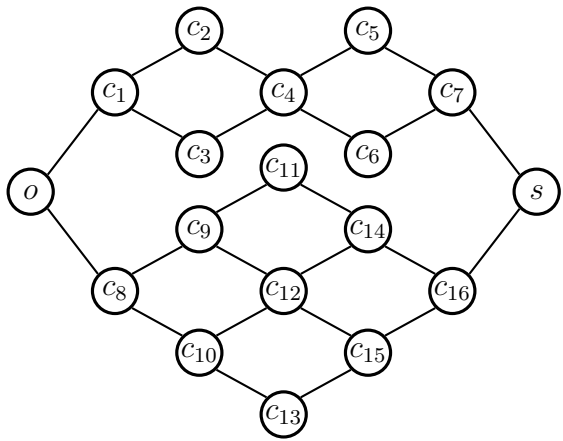

The second example is the 16-component network represented in \figref{fig:general_network_example}. The components have different topological importance: component $c_8$, $c_1$ and $c_4$, and the ones symmetric to them, can be considered as "bottlenecks", with respect to the other components.

\begin{figure}
        \centering
        \begin{subfigure}[b]{0.47\textwidth}
            \includegraphics[width=\textwidth]{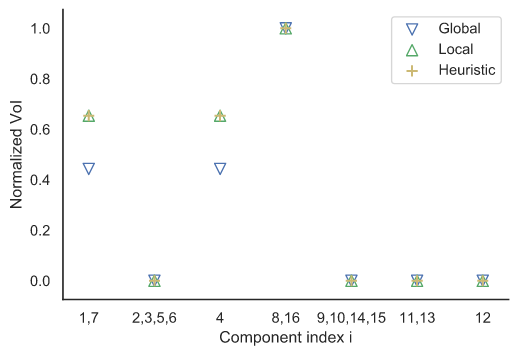}
            \caption{}
            \label{fig:VoI for network 1}
        \end{subfigure}
        \hfill
        \begin{subfigure}[b]{0.47\textwidth}
            \includegraphics[width=\textwidth]{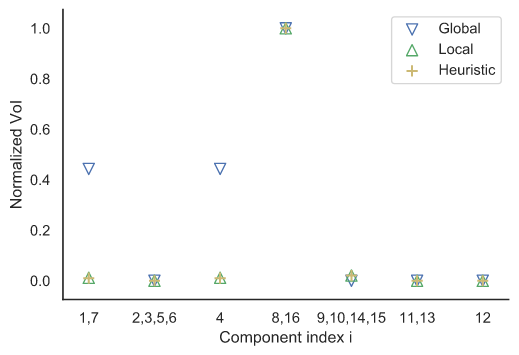}
            \caption{}
            \label{fig:VoI for network 2}
        \end{subfigure}
        \caption{Normalized VoI for the network in \figref{fig:general_network_example} with different maintenance cost, with and $C_F/C_{R,i}=10^{3}$ (a) and $C_F/C_{R,i}=10^{4}$ (b)}
        \label{fig:p001}
    \end{figure}

We assume that the marginal probabilities of failure is $p_i=0.01$ for every component $c_i$. For the global metric, we use $l^*(p_u)=p_u(1-p_u)$ as a loss function, and the corresponding VoI is shown in \figref{fig:VoI for network 1}. For the local metric, we assume that $C_F/C_{R,i}=10^3$ for every component $c_i$, so that the resulting optimal prior maintenance action is to repair no component. Under the local metric, $c_8$ and $c_{16}$ have the highest VoI, followed by $c_1$, $c_4$ and $c_7$. The heuristic approach gives the same result as the local metric, and the global metric is also consistent with this importance ordering. 

However, if $C_F/C_{R,i}$ is increased to $10^4$ for every component $c_i$, the new optimal prior maintenance action becomes to replace the symmetric bottlenecks $c_8$ and $c_{16}$. The VoI for the global and local metrics and the heuristic with this new assumption on costs is illustrated is \figref{fig:VoI for network 2}. The local VoI of inspecting $c_2$, $c_9$ and the components symmetric to them is now nil, because the cost for system failure is so (relatively) high, that the agent will not alter the prior action even if a silence is received on these components.

\begin{figure}
        \centering
        \includegraphics[width=.5\textwidth]{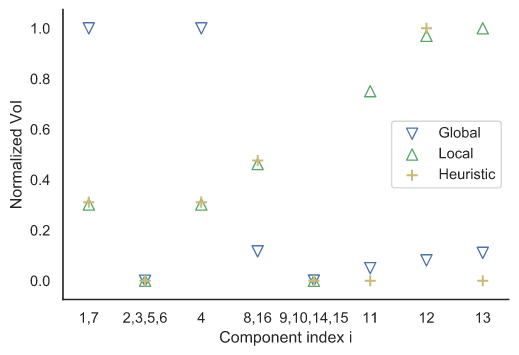}
        \caption{Normalized VoI for the network in \figref{fig:general_network_example} with $p_i=0.01,i\neq 11,12,13$, $p_{11}=0.5, p_{12}=0.4, p_{13}=0.3$ and $C_F/C_{R,i}=10^{3}$}
        \label{fig:VoI for network 3}
    \end{figure}

Depending on the setting, the bottleneck components may not always have the highest VoI. If $p_{11}=0.5, p_{12}=0.4, p_{13}=0.3$, $p_i=0.01,i\neq 11,12,13$ and $C_F/C_{R,i}=1000$ for every component $c_i$, the VoI is that shown in \figref{fig:VoI for network 3}. Now the optimal prior action is to replace $c_{12}$. The global metric prioritizes $c_1$, $c_4$ and $c_7$ for inspection, but the local metric prioritizes $c_{13}$, even though it is not the most vulnerable component (which is $c_{11}$). After $c_{13}$, the components with high VoI under local metric will be $c_{12}$ and $c_{11}$. Instead, the heuristic approach assigns the highest VoI to $c_{12}$. This is because when the inspection on $c_{11}$ or $c_{13}$ receives silence, the optimal action is to do nothing, but the heuristic approach forces the agent to at least execute the prior plan.

The third example is taken from \citep{Song2003BoundsReliability} and represents a two-line electrical substation with 12 components with 6 different functions as illustrated in \figref{fig:A two-transmission-line substation system}. In this example, we investigate how the correlation among component's state affects the VoI. If all the components are statistically dependent, complexity of computing the system failure probability may become intractable. \cite{Song2009} assumed conditional independence between component events given outcomes of a few random variables representing the source of common source effects and applied matrix-based method to compute system reliability. Following their framework, we assume interdependence among the components' states, but only for components with the same function. Concretely, the states of components $\text{DS}_1$, $\text{DS}_2$ and $\text{DS}_3$ are interdependent, while those of $\text{DS}_1$ and $\text{CB}_1$ are independent. The marginal failure probability of the components with function DS, CB or DB is $9.53\times10^{-3}$, and that of components with function PT or TB is $2.32\times10^{-3}$. For every component $c_i$, costs are defined by ratio $C_F/C_{R_i}=1000$. 

\begin{figure}
    \centering
    \begin{tikzpicture}[
        roundnode/.style={circle, draw=black!100, very thick, minimum width=9mm, inner sep=0pt},
        diamondnode/.style={diamond, draw=black!100, very thick, minimum width=5mm},
        squarenode/.style={rectangle, draw=black!100, very thick, minimum size=7mm},
        ]
        %Nodes
        \node[roundnode] (o) { $o$};
        \node[coordinate, right=0.8cm of o] (oright) {};
        \node[roundnode, above=1cm of oright] (n1) {\small $\text{DS}_1$};
        \node[roundnode, below=1cm of oright] (n2) {\small$\text{DS}_2$};
        \node[roundnode, right=0.8cm of oright] (n3) {\small$\text{DS}_3$};
        \node[coordinate, right=0.8cm of n1] (n1right) {};
        \node[coordinate, right=0.8cm of n2] (n2right) {};
        \node[roundnode, right=1cm of n1right] (n4) {\small$\text{CB}_1$};
        \node[roundnode, right=1cm of n2right] (n5) {\small$\text{CB}_2$};
        \node[roundnode, right=1cm of n4] (n6) { \small$\text{PT}_1$};
        \node[roundnode, right=1cm of n5] (n7) { \small$\text{PT}_2$};
        \node[roundnode, right=1cm of n6] (n8) {\small$\text{DB}_1$};
        \node[roundnode, right=1cm of n7] (n9) {\small$\text{DB}_2$};
        \node[coordinate, right=0.8cm of n8] (n8right) {};
        \node[coordinate, right=0.8cm of n9] (n9right) {};
        \node[roundnode, right=0.8cm of n8right] (n11) {\small$\text{FB}_1$};
        \node[roundnode, right=0.8cm of n9right] (n12) {\small$\text{FB}_2$};
        \node[roundnode, below=1cm of n8right] (n10) {TB};
        \node[roundnode, right=1.6cm of n10] (s) { $s$};
        
        \draw[thick] (o.north east) -- (n1.south west);
        \draw[thick] (o.south east) -- (n2.north west);
        \draw[thick] (n1.east) -- (n4.west);
        \draw[thick] (n4.east) -- (n6.west);
        \draw[thick] (n6.east) -- (n8.west);
        \draw[thick] (n8.east) -- (n11.west);
        \draw[thick] (n2.east) -- (n5.west);
        \draw[thick] (n5.east) -- (n7.west);
        \draw[thick] (n7.east) -- (n9.west);
        \draw[thick] (n9.east) -- (n12.west);
        \draw[thick] (n3.north) -- (n1right);
        \draw[thick] (n3.south) -- (n2right);
        \draw[thick] (n10.north) -- (n8right);
        \draw[thick] (n10.south) -- (n9right);
        \draw[thick] (n11.south east) -- (s.north west);
        \draw[thick] (n12.north east) -- (s.south west);
        \node[text width=10.5cm,below=0.3cm of n7.south] {\small DS: Disconnect Switch, CB: Circuit Breaker, PT: Power Transformer, DB: Drawout Breaker, TB: Tie Breaker, FB: Feeder Breaker};
    \end{tikzpicture}
    \caption{Scheme of a two-transmission-line substation system}
    \label{fig:A two-transmission-line substation system}
\end{figure}

\begin{figure}[htbp]
        \centering
        \begin{subfigure}[b]{0.47\textwidth}
            \includegraphics[width=\textwidth]{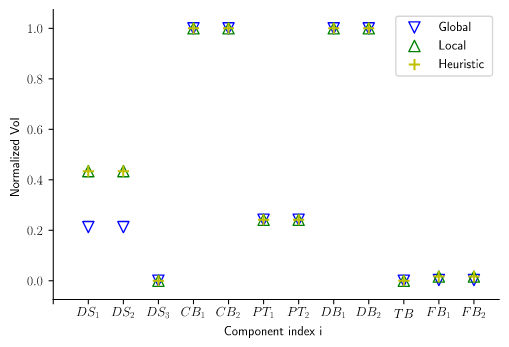}
            \caption{}
            \label{fig:corr0}
        \end{subfigure}
        \hfill
        \begin{subfigure}[b]{0.47\textwidth}
            \includegraphics[width=\textwidth]{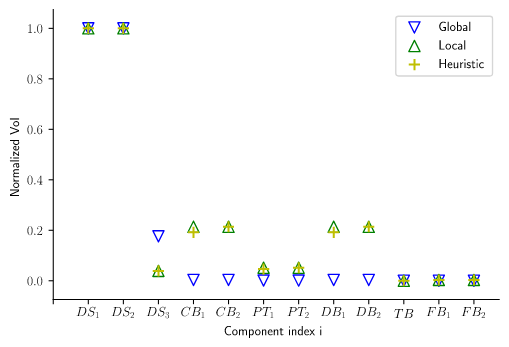}
            \caption{}
            \label{fig:corr02}
        \end{subfigure}
        \begin{subfigure}[b]{0.47\textwidth}
            \includegraphics[width=\textwidth]{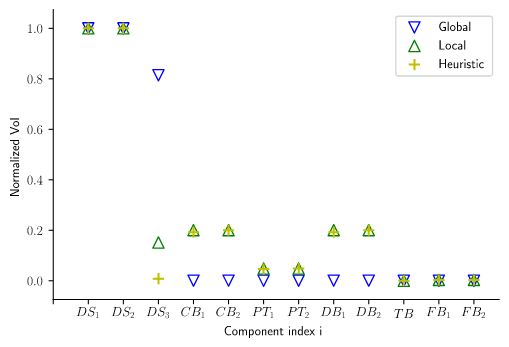}
            \caption{}
            \label{fig:corr06}
        \end{subfigure}
        \caption{Normalized VoI for the network in \figref{fig:A two-transmission-line substation system}, with correlation among the DS component of $\rho=0$ (a),  $\rho=0.4$ (b), $\rho=0.9$ (c) }
        \label{fig:VoI for network corr}
    \end{figure}

When all the components are independent, the prior action is to do nothing, and the optimal posterior action is to replace the inspected component after an alarm, except for $\text{DS}_3$ and $\text{TB}$. Thus, the local metric and heuristic give identical results. Though $\text{CB}$ and $\text{DB}$ have relatively higher failure probability than other components, the cost reduction by replacing the damaged components $\text{CB}$ or $\text{DB}$ is significantly higher than others. This is why these components have highest VoI according to the local metric and the heuristic, as shown in \figref{fig:corr0}. For the global metric, with loss function $l^*(p_u)=p_u(1-p_u)$, the posterior system failure probability given an alarm from components $\text{CB}$ or $\text{DB}$ is the highest, and the probability given a silence from those components is the lowest, i.e. posterior intervals $I(\text{CB})=I(\text{DB})$ contain the corresponding intervals of all the others, thus those components have the highest VoI, according to the global metric.

When the correlation among states in $\text{DS}$ components grows, while other groups remain independent (and the marginal probability remains the same), the VoI favors the group of correlated components. The prior action becomes replacing $\text{DS}_1$ or $\text{DS}_2$ when the correlation coefficient $\rho$ is above $0.4$. The optimal action is shown in \tabref{tab:Optimal posterior action for the substation system example 2}. Components $\text{DS}_1$ or $\text{DS}_2$ should be kept functioning, depending on which link set the inspected component is in. One exception is $\text{DS}_3$, which has different VoI for the local metric and the heuristic. As shown in \figref{fig:corr02}, when the correlation coefficient $\rho$ for the states of components $\text{DS}$ increases, inspecting one of them reveals additional information about the other two, making the VoI higher than that for other independent components. When $\rho$ is close to one, $DS_1$ and $DS_2$ act like one bottleneck component, which dominates the VoI as shown in  \figref{fig:corr06}.

\begin{table}[htbp]
    \centering
    \begin{tabular}{ccc}
        \toprule
        \diagbox[width=10em]{Component}{Insp. outcome}&
        Silence ($y_i=1$) & Alarm ($y_i=0$)\\
        \midrule
        $\text{DS}_1$ & $\emptyset$ & $\text{DS}_1$ \\
        $\text{DS}_2$ & $\emptyset$ & $\text{DS}_2$ \\
        $\text{DS}_3$ & $\emptyset$ & $\text{DS}_3$ \\
        $\text{CB}_1$ & $\text{DS}_1$ & $\text{DS}_1$, $\text{CB}_1$\\
        $\text{CB}_2$ & $\text{DS}_2$ & $\text{DS}_2 $, $\text{CB}_2$ \\
        $\text{PT}_1$ & $\text{DS}_1$ & $\text{DS}_1$, $\text{PT}_1$ \\
        $\text{PT}_2$ & $\text{DS}_2$ & $\text{DS}_2$, $\text{PT}_2$ \\
        $\text{DB}_1$ & $\text{DS}_1$ & $\text{DS}_1$, $\text{DB}_1$ \\
        $\text{DB}_2$ & $\text{DS}_2$ & $\text{DS}_2$, $\text{DB}_2$ \\
        $\text{TB}$ & $\text{DS}_1$ & $\text{DS}_1$ \\
        $\text{FB}_1$ & $\text{DS}_1$ & $\text{DS}_1$, $\text{FB}_1$ \\
        $\text{FB}_2$ & $\text{DS}_2$ & $\text{DS}_2$, $\text{FB}_2$ \\
        \bottomrule
    \end{tabular}
    \caption{Optimal posterior action for the network in \figref{fig:A two-transmission-line substation system} when $\rho=0.4$}
    \label{tab:Optimal posterior action for the substation system example 2}
\end{table}

\section{Conclusion}

We have derived metrics based on the VoI to assign priorities among component inspections in network systems. The VoI analysis can be applied to any setting, but its computational complexity depends on the complexity of the ingredients defining the problem. We have restricted the attention to binary components, binary inspection outcomes in a binary system. In this setting, we have introduced two metrics, the local and the global ones, that assume different sets of available actions and, thus, different loss functions. The problems modeled by the global metric do not form a sub-class of that of the problems modeled by the local one. We have proven general rules of identifying what components have higher importance for those metrics, in some simple systems.
The evaluation of the global metric is generally less complex than that of the local one, as so is the underlying optimization of the maintenance actions.  The selection of the appropriate metric should be based on the actual set of actions available. However, when only limited computational resources are available, the selection of a simpler metric or of the heuristic can be needed or appropriate.

We have proposed a heuristic approach to approximate the local metric, by simplifying the corresponding optimization of maintenance actions. There is no guarantee that the heuristic captures the exact local metric, and its performance has been illustrated in some examples. The VoI assessed by the heuristic is surely non-negative, and no higher than that of the original local metric, however the ranking can be arbitrarily different.

The metrics can be extended to more complex settings. The distinction between local and global metric is summarized in \figref{fig:figure1}: if the actions directly affect the system state $u$, so that a concave function can be defined on the domain of that variable, then the problem refers to the global metric. Otherwise, if actions affect the joint state $s$ of all components, it refers to the local metric. This distinction can be extended to the case of multiple values (more than binary) for the state of the components and of the system, and for inspection outcomes. However, some concepts are defined only for problems in small dimensions, e.g. the posterior intervals in the global metric are defined only for binary inspection outcomes in a binary system.

\section*{Acknowledgement}
The first and second authors acknowledge the support of NSF project CMMI 1653716, titled “CAREER: Infrastructure Management under Model Uncertainty: Adaptive Sequential Learning and Decision Making.” The second author acknowledges the support of Visiting Faculty Fellows Program of the Wilton E. Scott Institute for Energy Innovation at Carnegie Mellon University.

\appendix
\section{Importance Measures}
\label{app:Importance_Measures}

Similar to the Birnbaum's measure, the Criticality IM \citep{espiritu2007component}, evaluates the importance of $c_i$ with the approximated conditional component failure probability given that the system has failed:
\begin{equation}
    \text{CRT}(i)= (p_{\omega|\mathit{y}_i=0} - p_{\omega|\mathit{y}_i=1})\frac{p_i}{p_\pi} \propto \text{BM}(i)\cdot p_i
\end{equation}

Some IMs put emphasis on the topology structure of the network. Based on the cut sets of the network, \citep{fussell1975hand} evaluates the importance of $c_i$ by the number of cut sets it belongs to and the accumulated appearance probability of such cut sets. 

To use IMs as utility-based applications, the risk achievement worth (RAW) and the risk reduction worth (RRW) are developed. RAW evaluates the component with the contributions of maintaining a certain level of reliability of the component to the system reliability, i.e. for component $c_j$, its importance can be measured as:
\begin{equation}
    \text{RAW}(i) = \frac{1-p_{\omega|\mathit{y}_i=1}}{p_\pi}
\end{equation}
So, between two components $c_i$ and $c_j$, $\text{RAW}(i) \geq \text{RAW}(j)\Leftrightarrow p_{\omega|y_i=1}\leq p_{\omega|y_j=1}$. RRW evaluates a component by the decrease of system failure risk given that the component is intact:
\begin{equation}
    \text{RRW}(j) = \frac{p_\pi}{1-p_{\omega|\mathit{y}_i=0}}
\end{equation}
So $\text{RRW}(i) \geq \text{RRW}(j)\Leftrightarrow p_{\omega|y_i=0}\geq p_{\omega|y_j=0}$.

\section{Nested posterior intervals in the global metric}
\label{app:nested}

To prove the lemma in \secref{sub:priority_for_nested}, we now write $p_{\omega|y_a=b}$ as $x_{a,b}$ for simplicity. We assume that $I_i\supseteq I_j$, we have that $0\leq x_{i,1}\leq x_{j,1}\leq x_{j,0}\leq x_{i,0}\leq 1$. Because of the law of expectation, we have:
\begin{equation}
    p_\pi=p_1 x_{i,1} + (1-p_1) x_{i,0} = p_2 x_{j,1} + (1-p_2) x_{j,0}
\end{equation}
We prove that:
\begin{equation}
    L^\text{G}_\omega(1)=p_1 l(x_{i,1}) + (1-p_1) l(x_{i,0}) \leq p_2 l(x_{j,1}) + (1-p_2) l(x_{j,0})=L^\text{G}_\omega(2)
\end{equation}
Because $x_{j,1}=\frac{x_{i,0}-x_{j,1}}{x_{i,0}-x_{i,1}}x_{i,1} + \frac{x_{j,1}-x_{i,1}}{x_{i,0}-x_{i,1}}x_{i,0}$ and $x_{j,0}=\frac{x_{i,0}-x_{j,0}}{x_{i,0}-x_{i,1}}x_{i,1} + \frac{x_{j,0}-x_{j,1}}{x_{i,0}-x_{i,1}}x_{i,0}$, and $l$ is a concave function, we have:
\begin{equation}
    \begin{aligned}
        p_2 l(x_{j,1}) + (1-p_2) l(x_{j,0}) 
        \geq & p_2 [\frac{x_{i,0}-x_{j,1}}{x_{i,0}-x_{i,1}} l(x_{i,1}) + \frac{x_{j,1}-x_{i,1}}{x_{i,0}-x_{i,1}} l(x_{i,0})] \\
        &+ (1-p_2)[\frac{x_{i,0}-x_{j,0}}{x_{i,0}-x_{i,1}} l(x_{i,1}) + \frac{x_{j,0}-x_{j,1}}{x_{i,0}-x_{i,1}}l(x_{i,0})] \\
        =& p_1 l(x_{i,1}) + (1-p_1) l(x_{i,0})
    \end{aligned}
\end{equation}

%\end{linenumbers}

\bibliographystyle{elsarticle-harv}
\bibliography{reference.bib}

\end{document}